# Stochastic impulse control of non-smooth dynamics with partial observation and execution delay: application to an environmental restoration problem

Short title: Stochastic impulse control under partial observation with execution delay

Hidekazu Yoshioka[1, 2, *], Yuta Yaegashi[3]

[1] Assistant Professor, Graduate School of Natural Science and Technology, Shimane University, Nishikawatsu-cho 1060, Matsue, 690-8504, Japan

[2] Center Member, Fisheries Ecosystem Project Center, Shimane University, Nishikawatsu-cho 1060, Matsue, 690-8504, Japan

[3] Independent Researcher, Dr. of Agr., 10-12-403, Maeda-cho, Niihama, 792-0007, Japan

* Corresponding author
  E-mail: yoshih@life.shimane-u.ac.jp

**Abstract**

Non-smooth dynamics driven by stochastic disturbance arise in a wide variety of engineering problems. Impulsive interventions are often employed to control stochastic systems; however, the modeling and analysis subject to execution delay have been less explored. In addition, continuously receiving information of the dynamics is not always possible. In this paper, with an application to an environmental restoration problem, a continuous-time stochastic impulse control problem subject to execution delay under discrete and random observations is newly formulated and analyzed. The dynamics have a non-smooth coefficient modulated by a Markov chain, and eventually attain an undesirable state like a depletion due to the non-smoothness. The goal of the control problem is to find the most cost-efficient policy that can prevent the dynamics from attaining the undesirable state. We demonstrate that finding the optimal policy reduces to solving a non-standard system of degenerate elliptic equations, the Hamilton–Jacobi–Bellman equation, which is rigorously and analytically verified in a simplified case. The associated Fokker–Planck equation for the controlled dynamics is derived and solved explicitly as well. The model is finally applied to numerical computation of a recent river environmental restoration problem. The Hamilton–Jacobi–Bellman and Fokker–Planck equations are successfully computed, and the optimal policy and the probability density functions are numerically obtained. The impacts of execution delay are discussed to deeper analyze the model.

## 1. Introduction

We consider a stochastic control problem of a non-smooth dynamical system with partial observation and delayed execution: a new problem related to an environmental issue. In this section, we firstly review the stochastic control and then explain our objectives and contributions.

Optimal control theory is a pillar of modern mathematical sciences dealing with management of dynamical systems in engineering problems [1]. Stochastic control is a branch of the optimal control specialized for problems with noise-driven dynamics [2]. Management problems not limited to but include those of energy and resources [3, 4], population [5, 6], environment [7, 8], finance and economics [9, 10], and planning [11, 12] have been analyzed as stochastic control problems.

It is widely accepted that continuous-time stochastic control models based on stochastic differential equations (SDEs) [13, 14] serve as an efficient tool for approaching many scientific and engineering problems. In a standard modeling strategy of the stochastic control, some SDEs to be controlled are formulated and a performance index (or called objective function) to be optimized by controlling the dynamics is specified. Then, applying a dynamic programming principle [2] reduces the control problem to a degenerate parabolic or elliptic equation, called the Hamilton–Jacobi–Bellman equation (HJBE), from which the optimal control would be found. Solving the HJBE is in general carried out either analytically [15-17] or numerically [18-20].

The conventional framework assumes that the decision-maker controlling the target dynamics can receive complete information of the dynamics. However, this assumption is often violated in real problems especially in those related to environment and ecology. Continuously observing environmental and biological dynamics is a difficult task, and scheduled discrete observations of environmental and biological dynamics are more reasonable and feasible alternatives [21-23]. Under such situations, the decision-maker must make decisions based only on partial information, and the HJBE becomes more complicated because of using strictly smaller filtrations [24, 25]. In some models, even the assumption of the scheduled observations is relaxed and only random observation processes are assumed [26]. This partial observation framework leads to a control problem subject to discrete and random observations that has been used in some financial models under restricted observation chances [27-29].

Another difference between the ideal and real problems is the execution delay. Interventions that have much shorter time-scales than those of the target dynamics are reasonably considered to be impulsive [2]. Most of the impulse control models assume that the interventions are executed immediately at making the decisions [30, 31], while there would be execution delays in real problems due to technical reasons like communication and implementation delays. So far, problems with the execution delay have been studied theoretically from classical [32, 33] and viscosity viewpoints [34, 35] with dynamic programming [36]. Semi-analytical numerical schemes for approximating solutions to these problems have been developed in Perera and Long [37]. However, most of the existing models focus on financial and economic problems, but not those related to environment and ecology.

The objective of this paper is to formulate a stochastic impulse control problem based on discrete and random observation subject to execution delay, with a focus on a river environmental restoration problem by wisely using earth and soils [23, 38]. These previous studies considered coupled algae-sediment dynamics, while we solely focus on the sediment dynamics to build a simpler model such that we can handle it analytically under certain assumptions. The mathematical problem analyzed in this paper has not been explored well to the best of the authors' knowledge. The core of the problem is replenishing the sediment storage in a dam-downstream environment, which follows non-smooth dynamics modulated by a Markov chain representing river flows. We assume that the decision-maker can replenish the sediment impulsively with execution delay based on discrete and random observations. Namely, we only allow for interventions at the observation times, effectively making the problem a discrete time decision-making problem. We formulate the model under a more general setting to cover wider range of replenishment problems where the occurrence of an undesirable state such as a depletion state should be avoided.

We show that the problem above can be formulated based on a dynamic programming principle reducing a control problem to an HJBE as a system of degenerate elliptic equations. We demonstrate that the equation has a closed-form solution with some coefficients determined from algebraic equations under a simplification, and that a threshold type control is optimal by a verification argument. We also analyze the probability density function (PDF) of the controlled dynamics. Behavior of the PDFs of impulsively controlled stochastic systems has not been paid attention except for several recent studies [39, 40]. We heuristically derive the Fokker–Planck equation (FPE) governing the PDF of the controlled dynamics formulated based on the stochastic resetting [41] and controlled jump processes [42], and solve it analytically under the simplified case. The FPE is then verified with a Monte-Carlo simulation. Problems without the simplification are numerically handled using high-resolution finite difference schemes equipped with Weighted Essentially Non-Oscillatory (WENO) reconstructions [43, 44]. We show that they successfully approximate the exact solutions to the two equations under the simplification. The optimal sediment replenishment policy is finally analyzed numerically and the impacts of execution delay and model aggregation are discussed.

In summary, we contribute to formulation, analysis, computation, and application of a stochastic impulse control problem with execution delay and random discrete observation. The rest of this paper is organized as follows. The control problem is formulated and the optimality is presented in Section 2. A tractable case is studied in Section 3. The presented model is applied to numerical computation of a real problem in Section 4. Summary of this paper and future perspectives of our research are presented in Section 5. **Appendices** contain technical proofs of the analysis results and an explanation of the numerical scheme for the FPE.

## 2. Mathematical model
### 2.1 Dynamics without interventions

We consider a system governing a continuous-time variable $X = (X_t)_{t \geq 0}$ driven by a continuous-time Markov chain $\alpha = (\alpha_t)_{t \geq 0}$. Here, $t$ represents the time. The total number of the regimes of $\alpha$ is $I + 1 \in \mathbb{N}$ and is numbered from $i = 0$ to $i = I$. As in the standard setting of Markov chains [14], the switching matrix of $\alpha$ is set as $v = [v_{i,j}]_{i,j \in M}$ with non-negative entries; set $M = \{i\}_{0 \leq i \leq I}$. We assume that $v$ is irreducible and all the regimes are transient.

Assume that $X$ represents a system state such that a larger value corresponds to a more desirable state and vice versa, and that the most undesirable state $X_t = 0$ persists once it is attained unless there is some replenishment. Set the range of $X$ as $D = [0,1]$ with a normalization. The degradation speed of the state is represented by $S : M \times D \to \mathbb{R}_+$ such that $S$ is non-negative, bounded in $M \times D$, $S(\cdot, 0) = 0$, and satisfies the one-sided Lipschitz continuity for each $i \in M$:

$$-(x_1 - x_2)[S(i, x_1) - S(i, x_2)] \leq L(x_1 - x_2)^2 \quad \text{for} \quad x_1, x_2 \in D, \tag{1}$$

where $L \geq 0$ is a constant. This is a standard assumption to well-pose non-smooth dynamics [45]. For example, $S(i, x) = S_i \chi_{\{x > 0\}}$ satisfies (1), where $S_i > 0$ is a constant and $\chi_A$ is the indicator function for the set $A$ (e.g., $\chi_{\{x > 0\}} = 1$ if $x > 0$ and $\chi_{\{x = 0\}} = 0$ otherwise).

Set $X_0 = x \in D$ and the governing SDE of $X$ without interventions as

$$\mathrm{d}X_t = -S(\alpha_{t-}, X_t)\mathrm{d}t \quad \text{for} \quad t > 0. \tag{2}$$

This equation represents that the state degrades with the speed $S$. Its continuous solution can be formally expressed as

$$X_t = \max\left\{0, x - \int_0^t S(\alpha_{s-}, X_s)\mathrm{d}s\right\} \quad \text{for} \quad t \geq 0. \tag{3}$$

Because we are interested in a control problem to avoid $X_t = 0$, we assume that the state eventually becomes undesirable:

$$\int_0^t S(\alpha_{s-}, X_s)\mathrm{d}s \to +\infty \quad \text{almost surely (a.s.) as} \quad t \to +\infty \tag{4}$$

and that (3) is the unique continuous path-wise solution to (2). These are satisfied for example if $S(i, x) = S_i \chi_{\{x > 0\}}$ and $S_{i_0} > 0$ with some $i_0 \in M$ [23, 46]. The solution is relevant for considering the replenishment problem because of its non-increasing nature.

## 2.2 Dynamics with interventions

We assume that the decision-maker can replenish the state $X$ by impulsive interventions. The observation times are represented by an increasing sequence $\tau = \{\tau_k\}_{k \in \mathbb{N}}$ and $\tau_0 = 0$. The exponential distribution with a parameter $p$ is expressed as $\mathrm{Exp}_p$ and a random variable $r$ following $\mathrm{Exp}_p$ as $r \sim \mathrm{Exp}_p$.

Assume that the interventions are determined at each observation time with the following rule: do nothing ($\theta_k = 0$) or replenish ($\theta_k = 1$). If $\theta_k = 0$, then the dynamics are not impulsively controlled at $\tau_k$ and $\tau_{k+1} - \tau_k \sim \mathrm{Exp}_\lambda$ with $\lambda > 0$. On the other hand, if $\theta_i = 1$, then the state is impulsively moved to the most desirable state $X_{\tau_k + \omega_k} = 1$ at a future time $\tau_k + \omega_k$, where $\omega_k \sim \mathrm{Exp}_\mu$ with $\mu(> \lambda)$ and $\tau_{k+1} - (\tau_k + \omega_k) \sim \mathrm{Exp}_\lambda$. Set $\Delta_k = 0$ if $\theta_k = 0$ and $\Delta_k = \omega_k$ if $\theta_k = 1$, and $\theta_0 = 0$. Without significant loss of generality, we assume that each waiting time, each delay time, and the Markov chain $\alpha$ are mutually independent with each other. We are assuming that each delay ($\omega_k$) and waiting time ($\tau_{k+1} - (\tau_k + \Delta_k)$) follow some exponential distributions and that the delay times are likely to be shorter than the waiting times. The exponential delay assumptions have been utilized to construct simpler models with delayed execution [32, 33]. Specifically, the present model becomes tractable owing to utilizing the exponential distributions (See, Section 3).

Consequently, the decision-maker can discretely and randomly observe the dynamics, and the state is impulsively updated with a random delay (**Figure 1**).

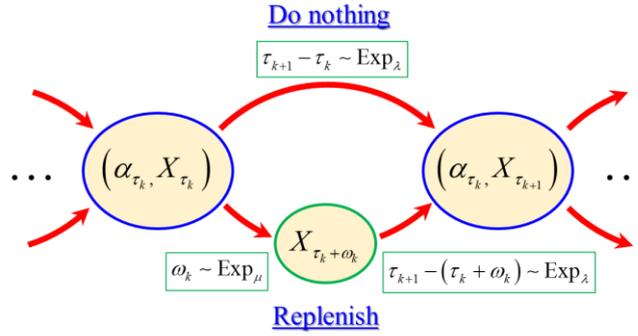

**Figure 1**. Sketch of the observations and interventions.

The decision-maker can observe each $\left(\tau_j, \alpha_{\tau_j}, X_{\tau_j}\right)$ and $\left(\tau_j + \Delta_j, \alpha_{\tau_j + \Delta_j}, X_{\tau_j + \Delta_j}\right)$, suggesting to set a natural filtration generated by the available information as $\mathcal{F} = (\mathcal{F}_t)_{t \geq 0} = (\mathcal{F}_{1,t} \vee \mathcal{F}_{2,t})_{t \geq 0}$, where

$$\mathcal{F}_{1,t} = \sigma\left\{\left(\tau_j, \alpha_{\tau_j}, X_{\tau_j}\right)_{0 \leq j \leq k}, k = \sup\{j : \tau_j \leq t\}\right\} \tag{5}$$

and

$$\mathcal{F}_{2,t} = \sigma\left\{\left(\tau_j + \Delta_j, \alpha_{\tau_j + \Delta_j}, X_{\tau_j + \Delta_j}\right)_{0 \leq j \leq k}, k = \sup\{j : \tau_j + \Delta_j \leq t\}\right\}. \tag{6}$$

In summary, the dynamics of $X$ subject to the interventions are formulated as

$$dX_t = -S(\alpha_{t-}, X_t) dt \quad \text{for} \quad t > 0, \ t \neq \tau_k + \Delta_k \ (k \in \mathbb{N}) \tag{7}$$

with

$$X_{t+} = \begin{cases} 1 & (\theta_k = 1) \\ X_t & (\theta_k = 0) \end{cases} \quad \text{for} \quad t = \tau_k + \Delta_k \ (k \in \mathbb{N}). \tag{8}$$

The amount of replenishment at each $\tau_k + \Delta_k$ is

$$\eta_k = \begin{cases} 1 - X_{\tau_k + \Delta_k} & (\theta_k = 1) \\ 0 & (\theta_k = 0) \end{cases}. \tag{9}$$

A set of admissible control $C$ contains continuous-time processes $\bar{\eta} = (\bar{\eta}_t)_{t \geq 0}$ such that it equals $\eta_k$ at each $\tau_k + \Delta_k$ and 0 otherwise, and $\eta_k$ is adapted to $\mathcal{F}_{\tau_k}$. The last condition is due to that the replenishment may be carried out at $t = \tau_k + \Delta_k$ but its decision at $t = \tau_k$.

## 2.3 Performance index and value function

A performance index is an index to be minimized with respect to $\bar{\eta} \in C$. The undesirable state $X_t = 0$ can be avoided by interventions, but such activities are costly. We consider an infinite-horizon setting to formulate a simpler problem. The conditional expectation with respect to $(\alpha_0, X_0) = (i, x)$ is denoted as $\mathbb{E}^{i,x}$. The performance index to be minimized is

$$\phi(i, x, \bar{\eta}) = \mathbb{E}^{i,x}\left[\int_0^\infty e^{-\delta s}\chi_{\{X_s = 0\}} \mathrm{d}s + \sum_{k \geq 1} e^{-\delta(\tau_k + \Delta_k)}\left(c\eta_k + d\chi_{\{\eta_k > 0\}}\right)\right], \quad (10)$$

where $\delta > 0$ is the discount rate representing a myopic decision-maker by a larger value, $c > 0$ is the coefficient of proportional cost, and $d > 0$ is the fixed cost. The first term in the right-hand side of (10) penalizes the undesirable state, and the second term represents the intervention costs paid at the time of replenishment. No cost is incurred if there is no replenishment ($\eta_k = 0$).

Finally, the value function $\Phi : M \times D \to \mathbb{R}$ is the minimized $\phi$:

$$\Phi(i, x) = \inf_{\bar{\eta} \in C} \phi(i, x, \bar{\eta}). \quad (11)$$

The goal of the control problem is to find a minimizer $\eta^*$, an optimal control, of (11). Notice that $\Phi$ is well-defined since clearly $\phi \geq 0$ and the right-hand side of (10) is bounded because of $\delta > 0$ and the compactness of $D$. We can check this by substituting the null control ($\eta_k = 0$, $k \geq 1$) to (11).

## 2.4 HJBE

We present the HJBE governing $\Phi$. It is justified in the next section. Based on the heuristic arguments of the discrete and random observations [26] and execution delay [33], our HJBE is

$$\delta \Phi_i + S(i, x)\chi_{\{x > 0\}} \frac{\mathrm{d}\Phi_i}{\mathrm{d}x} + \sum_{j \in M, j \neq i} v_{i,j}(\Phi_i - \Phi_j) + \lambda\left(\Phi_i - \min\{\Phi_i, \hat{\Phi}_i\}\right) = \chi_{\{x = 0\}}, \quad (i, x) \in M \times D \quad (12)$$

and

$$\hat{\Phi}_i = \mathbb{E}^{i,x}\left[\int_0^\omega e^{-\delta s}\chi_{\{X_s = 0\}} \mathrm{d}s + e^{-\delta\omega}(c(1 - X_\omega) + d) + e^{-\delta\omega}\Phi_{\alpha_\omega}(1)\right], \quad (i, x) \in M \times D \quad (13)$$

with the notations $\Phi_i(x) = \Phi(i, x)$, $\hat{\Phi}_i(x) = \hat{\Phi}(i, x)$, and $\omega \sim \mathrm{Exp}_\mu$. No explicit boundary condition is necessary at $x = 0, 1$ since the characteristics of the process $X$ are vanishing or inward at the boundaries. Note that the right-hand side of (12), as well as the drift in (12), are discontinuous at only on one point $x = 0$, and thus that a mathematical difficulty, if it exists, is hidden here. Fortunately, in the present case, looking for a smooth solution should fix the value of its derivative at 0 and could help simplify our computations.

Heuristically, given an observation $(\alpha_{\tau_k}, X_{\tau_k})$ at $\tau_k$, assuming a Markov control of the form $\eta^* = \eta^*(\alpha_{\tau_k}, X_{\tau_k})$, we guess the optimal control

$$\eta^*(\alpha_{\tau_k}, X_{\tau_k}) = \begin{cases} 0 & (\Phi(\alpha_{\tau_k}, X_{\tau_k}) < \hat{\Phi}(\alpha_{\tau_k}, X_{\tau_k})) \\ 1 - X_{\tau_k + \omega_k} & (\Phi(\alpha_{\tau_k}, X_{\tau_k}) \geq \hat{\Phi}(\alpha_{\tau_k}, X_{\tau_k})) \end{cases}, \quad (14)$$

which is justified both mathematically and numerically in this paper. In this view, finding an optimal control is achieved by solving (12).

By $\omega \sim \mathrm{Exp}_\mu$, (13) is rewritten as

$$\begin{aligned}
\hat{\Phi}_i &= \mathbb{E}^{i,x}\left[\int_0^\infty \mu e^{-\mu t}\left\{\int_0^t e^{-\delta s}\chi_{\{X_s=0\}}\,\mathrm{d}s + e^{-\delta t}\left(c(1-X_t)+d\right) + e^{-\delta t}\Phi_{\alpha_t}(1)\right\}\mathrm{d}t\right] \\
&= \mathbb{E}^{i,x}\left[\int_0^\infty \int_0^t \mu e^{-\mu t}e^{-\delta s}\chi_{\{X_s=0\}}\,\mathrm{d}s\mathrm{d}t + \int_0^\infty \mu e^{-(\mu+\delta)t}\left\{\left(c(1-X_t)+d\right) + \Phi_{\alpha_t}(1)\right\}\mathrm{d}t\right] \\
&= \mathbb{E}^{i,x}\left[\int_0^\infty \int_s^\infty \mu e^{-\mu t}e^{-\delta s}\chi_{\{X_s=0\}}\mathrm{d}t\,\mathrm{d}s + \int_0^\infty \mu e^{-(\mu+\delta)t}\left\{\left(c(1-X_t)+d\right) + \Phi_{\alpha_t}(1)\right\}\mathrm{d}t\right] \\
&= \mathbb{E}^{i,x}\left[\int_0^\infty e^{-(\mu+\delta)t}\left\{\chi_{\{X_t=0\}} + \mu\left(c(1-X_t)+d\right) + \mu\Phi_{\alpha_t}(1)\right\}\mathrm{d}t\right]
\end{aligned} \quad (15)$$

where the order of integrations has been exchanged between the second and third lines. This representation of $\hat{\Phi}_i$ is more convenient from a computational viewpoint because invoking the Feynman-Kac formula [47] formally leads to the governing equation of $\hat{\Phi}$ given $\Phi$:

$$(\mu+\delta)\hat{\Phi}_i + S(i,x)\chi_{\{x>0\}}\frac{\mathrm{d}\hat{\Phi}_i}{\mathrm{d}x} + \sum_{j\in M, j\neq i} v_{i,j}(\hat{\Phi}_i - \hat{\Phi}_j) = \chi_{\{x=0\}} + \mu(c(1-x)+d) + \mu\Phi_i(1). \quad (16)$$

Again, no explicit boundary condition is necessary for $x = 0, 1$. Consequently, our HJBE is a coupled degenerate elliptic system containing (12) and (16). Stochastic control models with impulsive interventions with delayed executions sometimes associate HJBEs containing more than one degenerate elliptic equations [34]. However, to the best of the authors' knowledge, the above-presented HJBEs has not been presented so far.

## 3. Mathematical analysis
### 3.1 Optimal control under a simplified case

We show that the HJBE (12) gives an optimal control under a simplified condition. We analyze a single-regime case ($I = 0$) and $S(i,x) = S\chi_{\{x>0\}}$ with $S > 0$. The subscript $i$ is omitted here. The sub- and super-scripts representing the regimes are omitted here. A more realistic case is numerically analyzed in Section 4.

Under this setting, (12) simplifies to

$$\delta\Phi + S\chi_{\{x>0\}}\frac{\mathrm{d}\Phi}{\mathrm{d}x} + \lambda\left(\Phi - \min\{\Phi, \hat{\Phi}\}\right) = \chi_{\{x=0\}}. \quad (17)$$

By (15), $\hat{\Phi}$ is found as

$$\hat{\Phi}(x) = \alpha - c\beta x + \gamma e^{-\frac{\delta+\mu}{S}x} + \beta\Phi(1) \tag{18}$$

with constants

$$\alpha = \frac{\mu}{\delta+\mu}\left(c+d+\frac{cS}{\delta+\mu}\right), \quad \beta = \frac{\mu}{\delta+\mu}, \text{ and } \gamma = \frac{1}{\delta+\mu} - \frac{\mu cS}{(\delta+\mu)^2}. \tag{19}$$

A reasonable control policy would replenish the state if it is close to the undesirable state 0:

$$\eta_k^* = \begin{cases} 1 - X_{\tau_k+\omega_k} & (0 \leq X_{\tau_k} \leq \bar{x}) \\ 0 & (\bar{x} < X_{\tau_k} \leq 1) \end{cases} \tag{20}$$

with some threshold value $\bar{x} \in (0,1)$. We explore a smooth solution $\bar{\Phi} \in C^1(D)$ and show that it equals the value function $\Phi$ and further that this policy is optimal under certain assumption.

By (14), we infer

$$\delta\bar{\Phi} + S\chi_{\{x>0\}}\frac{d\bar{\Phi}}{dx} + \lambda(\bar{\Phi} - \hat{\Phi}) = \chi_{\{x=0\}} \text{ for } 0 \leq x \leq \bar{x} \text{ and } \delta\bar{\Phi} + S\frac{d\bar{\Phi}}{dx} = 0 \text{ for } \bar{x} < x \leq 1, \tag{21}$$

from which we get

$$\bar{\Phi}(x) = \begin{cases} De^{-\frac{\delta+\lambda}{S}x} + Ax + B + Ce^{-\frac{\delta+\mu}{S}x} & (0 \leq x \leq \bar{x}) \\ \bar{\Phi}(1)e^{\frac{\delta}{S}(1-x)} & (\bar{x} < x \leq 1) \end{cases} \tag{22}$$

with

$$A = \frac{\lambda\alpha}{\delta+\lambda} + \frac{\lambda cS\beta}{(\delta+\lambda)^2} + \frac{\lambda\beta\bar{\Phi}(1)}{\delta+\lambda}, \quad B = -\frac{\lambda c\beta}{\delta+\lambda}, \quad C = -\frac{\lambda}{\mu-\lambda}\gamma, \quad D = \frac{1}{\delta+\lambda}\left(1 - \frac{\lambda cS\beta}{\delta+\lambda} + \frac{\lambda(\delta+\mu)\gamma}{\mu-\lambda}\right). \tag{23}$$

There are the two unknowns $\bar{x}$ and $\bar{\Phi}(1)$, meaning that two equations are necessary to determine them. Assume $\bar{x} \in (0,1)$. We guess the smooth-pasting condition

$$\bar{\Phi}(\bar{x}-0) = \bar{\Phi}(\bar{x}+0) \text{ and } \frac{d\bar{\Phi}}{dx}(\bar{x}-0) = \frac{d\bar{\Phi}}{dx}(\bar{x}+0), \tag{24}$$

from which we try to determine the unknowns. They are lengthy and are not presented here (See, **Appendix B**), but can be directly obtained from (22). Although we have to solve the algebraic system (24) to get the two coefficients of the solution, the system can be solved using common numerical methods like standard Picard, Newton, or some trial and error methods with an arbitrary high accuracy.

Now, we present our first mathematical analysis result on the optimal control.

*Proposition 1*

*Assume that the algebraic equations (24) have a unique solution $(\bar{\Phi}(1), \bar{x}) \in \mathbb{R} \times (0,1)$. Then, we have the inequality $\Phi(x) \geq \bar{\Phi}(x)$ for all $x \in D$. Furthermore, the control of (20) is a minimizer of $\phi$ and gives $\Phi(x) = \bar{\Phi}(x)$ for all $x \in D$.*

Theoretically, the existence of the unique solution $(\bar{\Phi}(1), \bar{x})$ should be discussed to complete the statement of the existence of the simplified model. **Appendix B** addresses this issue analytically and shows that the threshold $\bar{x} \in (0,1)$ uniquely exists at least if the parameters $c, d, \mu^{-1}, \delta > 0$ are small; namely if the intervention costs are small, the delay intensity is small, and the discount rate is small. The last assumption means that the decision-maker is close to be Ergodic and thus controls the dynamics from a long-run viewpoint. As a byproduct of **Appendix B**, we get $\bar{\Phi}(1) > 0$. Then, (22) implies that the exact solution, the value function, is decreasing, smooth, and convex in $D$. See **Figure 2** plotting the value function, verifying the theoretical results. Especially, the conditions (14) and (20) are verified.

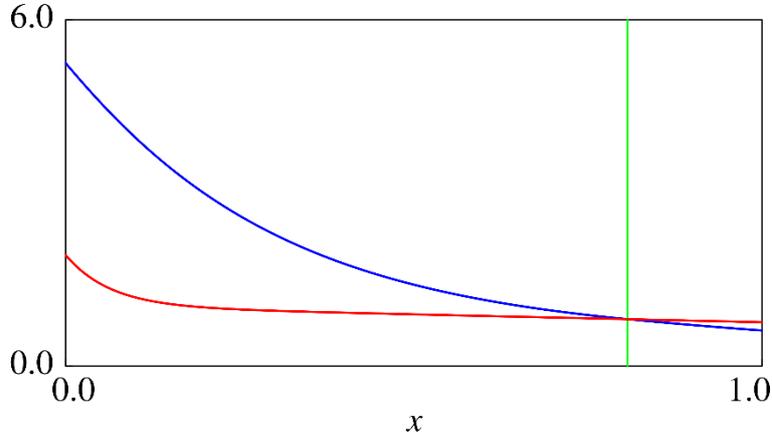

**Figure 2.** The functions $\bar{\Phi} = \bar{\Phi}(x)$ (Blue) and $\hat{\bar{\Phi}} = \hat{\bar{\Phi}}(x)$ (Red), and the corresponding threshold $\bar{x} = 0.807$ (Green). The parameter values have been chosen as follows: $\lambda = 1/7$, $\mu = 1$, $S = 0.07$, $c = 0.30$, $d = 0.20$, and $\delta = 0.1$.

### 3.2 Fokker–Planck equation and probability density function

An FPE equation is a governing equation of the PDF of state variables governed by a system of SDEs. For piecewise deterministic systems driven by Markov chains, the FPEs are systems of local linear hyperbolic equations [48, 49]. However, our FPE is a non-standard form because of the delayed interventions.

In this sub-section, we assume that the threshold value $\bar{x} \in (0,1)$ is given. There exist three issues to be considered in our problem. The first issue is that the state variable may be impulsively moved to a boundary point if it is smaller than $\bar{x}$ at the observation. The second issue is that this impulsive intervention is subject to the delay. The third issue is the non-smooth drift coefficient.

The first issue is resolved by invoking the stochastic resetting argument [41] handling the state variables moved to a point randomly. The second issue is then resolved by the jump rate argument [42] that can analytically characterize the jump rate from a state to another. In our case, the key fact is that the delay follows an exponential distribution. The last issue can be solved by considering an integral formulation. We formulate the FPE only heuristically as in Yaegashi et al. [40] and validate it using a Monte-Carlo method. Its derivation based on physical considerations is beyond the scope of this paper, and will be addressed elsewhere.

Assume a steady state condition and denote the PDF of the state $X_t = x$ as $p = p(x)$. There are two cases at each time instance; the state is during the delay duration waiting for the execution ($W$: Waiting case) or not ($N$: Non-waiting case). The PDFs of these cases are denoted as $p_W = p_W(x)$ and $p_N = p_N(x)$, respectively. We have $p_W + p_N = p$ and the mass conservation $\int_D p \, dx = 1$. Since the boundary points $x = 0, 1$ are not absorbing, the probability flux $F_a$, given by $-S\chi_{\{x>0\}} p_a$ inside the domain ($a$ represents $W$ or $N$), should vanish at $x = 0, 1$.

In this paper, the FPE of the controlled dynamics is formulated as

$$\frac{dF_N}{dx} + \lambda p_N = 0 \quad \text{for} \quad x \in [0, \bar{x}], \tag{25}$$

$$\frac{dF_N}{dx} - \mu c \delta_{\{x=1\}} = 0 \quad \text{for} \quad x \in (\bar{x}, 1] \quad \text{with} \quad c = \int_{[0,1]} p_W(x) \, dx, \tag{26}$$

$$\frac{dF_W}{dx} + \mu p_W - \lambda p_N = 0, \quad x \in [0, \bar{x}] \tag{27}$$

with

$$p_W(x) = 0 \quad \text{for} \quad x \in (\bar{x}, 1] \quad \text{and} \quad p_N(1) = 0, \tag{28}$$

where $\delta_{\{x=y\}}$ represents the Dirac Delta concentrated at $x = y$. Each equation is explained as follows. The first equation (25) means that the non-waiting case transits to the waiting case at the rate $\lambda$ if $x \leq \bar{x}$. The second equations (26) mean that the waiting case transits to the non-waiting case with $x = 1$ at the rate $\mu$ [41] only when the state just before was in the waiting case. The third equation (27) means that, when $x \leq \bar{x}$, the non-waiting case transits to the waiting case at the rate of $\lambda$ and the waiting case transits to the non-waiting case with $x = 1$ at the rate $\mu$. The last equations (28) mean that the waiting case with $x > \bar{x}$ does not occur due to the non-increasing nature of the state variable $X$, and the state variable immediately leave $x = 1$ by $S > 0$.

We give a remark on the boundary condition for the Non-waiting case that seems not to be trivial. After completing a replenishment event at time $t$, the state becomes $X_{t+} = 1$ in the Non-waiting case but the sediment storage immediately becomes smaller than 1 due to the negative drift. This physical observation implies that the sediment storage is smaller than 1 almost surely, suggesting to set $p_N(1) = 0$.

The equations (25) and (27) at $x = 0$ should be more clearly characterized due to the non-smoothness of the drift. Consider (25) in $[0, h]$ with $0 < h \ll \bar{x}$. An integral form of (25) is set as

$$\int_{[0,h]} \left( \frac{dF_N}{dx} + \lambda p_N \right) dx = 0. \tag{29}$$

Applying the classical divergence formula to (29) yields

$$0 = \int_{[0,h]} \left( \frac{dF_N}{dx} + \lambda p_N \right) dx = F_N(h) - F_N(0) + \lambda \int_{[0,h]} p_N dx = -Sp_N(h) + \lambda \int_{[0,h]} p_N dx. \tag{30}$$

By (30), and taking the limit $h \to +0$, we infer

$$p_N(+0) = \frac{\lambda}{S} \lim_{h \to +0} \int_{[0,h]} p_N dx, \tag{31}$$

implying that $p_N$ admits a Dirac mass at $x = 0$. If its proportional constant is denoted as $c_N \geq 0$, then

$$p_N(+0) = \frac{\lambda}{S} c_N. \tag{32}$$

Similarly, we infer from (30) that $p_W(0)$ is $c_W \delta_{\{x=0\}}$ with $c_W \geq 0$ satisfying

$$p_W(+0) = \frac{\mu}{S} c_W - \frac{\lambda}{S} c_N. \tag{33}$$

The appearance of the singular terms implies that the state $X_t = 0$ may persist until the next intervention.

The next proposition shows the existence of a PDF with $p_W, p_N \in C(0,1)$, showing in particular the positivity $c_N, c_W > 0$.

*Proposition 2*

*The FPE (25)-(28) admits $p_W, p_N \in C(0,1)$ with (32), (33), and $\int_D (p_N + p_W) dx = 1$. They are expressed as*

$$p_N(x) = \begin{cases} \mu \lambda^{-1} c e^{-\lambda S^{-1} \bar{x}} \delta_{\{x=0\}} & (x = 0) \\ \mu S^{-1} c e^{\lambda S^{-1}(x-\bar{x})} & (0 < x \leq \bar{x}) \\ \mu S^{-1} c (1 - \chi_{\{x=1\}}) & (\bar{x} < x \leq 1) \end{cases} \text{ and } p_W(x) = \begin{cases} Ec \delta_{\{x=0\}} & (x = 0) \\ Fc \left( e^{\lambda S^{-1}(x-\bar{x})} - e^{\mu S^{-1}(x-\bar{x})} \right) & (0 < x \leq \bar{x}) \\ 0 & (\bar{x} < x \leq 1) \end{cases} \tag{34}$$

*with*

$$E = \frac{\mu e^{-\lambda S^{-1} \bar{x}} - \lambda e^{-\mu S^{-1} \bar{x}}}{\mu - \lambda} > 0, \quad F = \frac{\lambda \mu}{S(\mu - \lambda)} > 0, \quad c = \frac{1}{1 + \mu \lambda^{-1} + \mu S^{-1}(1 - \bar{x})} > 0. \tag{35}$$

The above-presented exact solution is verified numerically. We use a Monte-Carlo method using $4 \times 10^6$ sample paths with a standard Forward-Euler method with the time increment of 0.0025. This is a purely statistical numerical method fundamentally different from the theoretical consideration above. The parameter values are set as $\lambda = 1/7$, $\mu = 1$, $S = 0.07$, and $\bar{x} = 0.807$. **Figure 3** compares the above-presented exact PDF and the computed one (histogram with the interval of 1/200) in $(0,1]$. The weights $c_W$ and $c_N$ at the origin $x = 0$ are 0.1253 and 0.02089 in the exact PDF, and 0.1254 and 0.02082 in the numerical computation, suggesting their good agreement. In addition, the computed profiles of the PDFs agree with each other from both qualitative and quantitative viewpoints. The computational results suggest that the FPE and the derived PDF are correct ones.

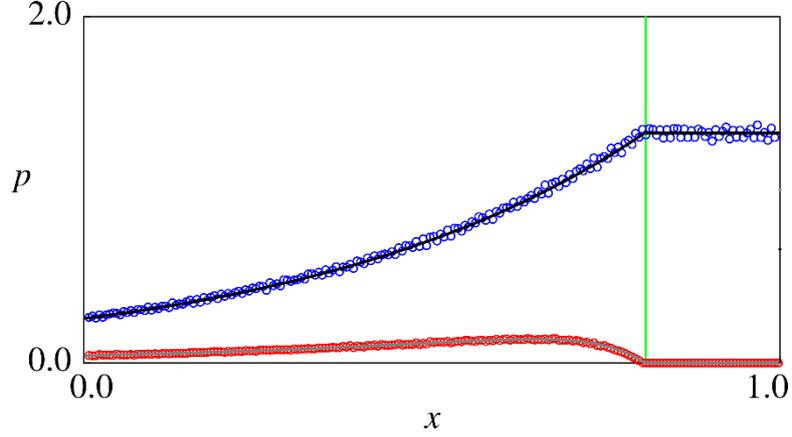

**Figure 3.** Comparison of the exact and numerical PDFs (Blue circles: non-waiting case with the Monte-Carlo method, Red circles: waiting case with the Monte-Carlo method, Black line: non-waiting case with the exact solution, Grey line: waiting case with the exact solution) in $(0,1]$. The green vertical line represents the threshold $\bar{x} = 0.807$.

### 3.3 The equations without the simplification

We conclude this section with remarks on the model without the above simplification. Firstly, the verification argument can be extended to the original HJBE (12). More specifically, assume the existence of a smooth solution $\Phi_i \in C^1(D)$, $i \in M$ with the threshold type control

$$\eta_k^* = \begin{cases} 1 - X_{\tau_k + \omega_k} & \left(0 \le X_{\tau_k} \le \bar{x}_{\alpha_{\tau_k}}\right) \\ 0 & \left(\bar{x}_{\alpha_{\tau_k}} < X_{\tau_k} \le 1\right) \end{cases} \quad (36)$$

with some $\bar{x}_i \in (0,1)$ ($i \in M$). Then, we can directly follow the argument of the **Proof of Proposition 1** (**Appendix A**) and obtain a similar verification result showing the optimality of (36). However, the existence of such a solution is not trivial because it is hopeless to get an exact solution explicitly as in the single-regime case. Later, we show numerically that the optimal policy is not always of the form (36), but maybe $\bar{x}_i = 0$ or 1, or even "Do nothing" is optimal for some $i \in M$.

Secondly, the FPE without the simplification can also be set heuristically. The difference between the problems with and without the simplification is that the former admits regime switching by the Markov chain $\alpha$, while it is absent in the latter. Based on this observation, the stationary PDFs in these cases are denoted as $p_{i,W} = p_W(i,x)$ and $p_{i,N} = p_N(i,x)$ for $(i,x) \in M \times D$, respectively. We should have the mass conservation $\sum_{i \in M} \int_D (p_{i,W} + p_{i,N}) dx = 1$. The probability flux is $F_{i,a} = -S(i,x) p_{i,a}$ ($a$ represents $W$ or $N$). The linear terms on the regime switching are then added to the FPEs [48, 49]. Consequently, for each $i \in M$, we obtain

$$\frac{dF_{i,N}}{dx} + \lambda p_{i,N} + \left(\sum_{j \in M, j \neq i} v_{i,j}\right) p_{i,N} - \sum_{j \in M, j \neq i} v_{j,i} p_{j,N} = 0 \quad \text{for} \quad x \in [0, \bar{x}], \quad (37)$$

$$\frac{dF_{i,N}}{dx} + \left(\sum_{j \in M, j \neq i} v_{i,j}\right) p_{i,N} - \sum_{j \in M, j \neq i} v_{j,i} p_{j,N} - \mu c_i \delta_{\{x=1\}} = 0 \quad \text{for} \quad x \in (\bar{x}, 1] \quad \text{with} \quad c_i = \int_{[0,1]} p_{i,W}(x) dx, \quad (38)$$

$$\frac{dF_{i,W}}{dx} + \mu p_{i,W} - \lambda p_{i,N} + \left(\sum_{j \in M, j \neq i} v_{i,j}\right) p_{i,W} - \sum_{j \in M, j \neq i} v_{j,i} p_{j,W} = 0, \quad x \in [0, \bar{x}] \quad (39)$$

with

$$p_{i,W}(x) = 0 \quad \text{for} \quad x \in (\bar{x}, 1] \quad \text{and} \quad p_{i,N}(1) = 0. \quad (40)$$

The singular nature at $x = 0$ would appear in this case as well, but the fact is more difficult to check analytically. We investigate this issue numerically along with behavior of the optimal control using identified parameter values and coefficients.

## 4. Application
### 4.1 Numerical schemes
#### 4.1.1 HJBE

The HJBE and the FPE are numerically discretized to compute the optimal control and the PDF of the controlled dynamics using realistic coefficients and parameter values. The numerical schemes used in this paper are finite difference and semi-Lagrangian schemes based on the WENO reconstructions. The WENO reconstructions have been widely employed to efficiently compute numerical solutions to linear and nonlinear degenerate elliptic equations and parabolic equations, owing to its unique nature of nonlinearly reconstructing solutions on computational grids so that numerical oscillations are minimized while the accuracy is maintained [50, 51].

The scheme for the HJBE is that of Carlini et al. [43] here applying a semi-Lagrangian discretization to $\delta \Psi_i + S(i,x) \chi_{\{x>0\}} \frac{d\Psi_i}{dx}$ ($\Psi$ represents $\Phi$ or $\hat{\Phi}$) and a classical forward Euler discretization to the other terms. The scheme itself is therefore not special, and is not explained in detail here. The HJBE is supplemented with artificial temporal terms as

$$\frac{\partial \Phi_i}{\partial t} + \delta \Phi_i + S(i,x) \chi_{\{x>0\}} \frac{\partial \Phi_i}{\partial x} + \sum_{j \in M, j \neq i} v_{i,j} (\Phi_i - \Phi_j) + \lambda \left(\Phi_i - \min\{\Phi_i, \hat{\Phi}_i\}\right) = \chi_{\{x=0\}} \quad (41)$$

and

$$\frac{\partial \hat{\Phi}_i}{\partial t} + (\mu + \delta) \hat{\Phi}_i + S(i,x) \chi_{\{x>0\}} \frac{d\hat{\Phi}_i}{dx} + \sum_{j \in M, j \neq i} v_{i,j} (\hat{\Phi}_i - \hat{\Phi}_j) = \chi_{\{x=0\}} + \mu (c(1-x) + d) + \mu \Phi_i(t,1), \quad (42)$$

where the solutions are now dependent on the (artificial) time $t \geq 0$ as $\Phi_i = \Phi_i(t,x)$ and $\hat{\Phi}_i = \hat{\Phi}_i(t,x)$. The original HJBE (12) is expected to be numerically recovered with a sufficiently large time $t \geq 0$. This pseudo-time discretization method has successfully been applied to degenerate elliptic equations [52]. The semi-Lagrangian nature of the scheme allows us to discretize the equation in a stable and efficient way [43]. The initial condition is set as $\Phi_i = \hat{\Phi}_i = 0$ in $M \times D$.

### 4.1.2 Fokker–Planck equation

The FPE is also computed by adding the artificial temporal partial differential terms:

$$\frac{\partial p_{i,N}}{\partial t} + \frac{\partial F_{i,N}}{\partial x} + \lambda p_{i,N} + \left( \sum_{j \in M, j \neq i} v_{i,j} \right) p_{i,N} - \sum_{j \in M, j \neq i} v_{j,i} p_{j,N} = 0 \text{ for } x \in [0, \bar{x}], \quad (43)$$

$$\frac{\partial p_{i,N}}{\partial t} + \frac{\partial F_{i,N}}{\partial x} + \left( \sum_{j \in M, j \neq i} v_{i,j} \right) p_{i,N} - \sum_{j \in M, j \neq i} v_{j,i} p_{j,N} - \mu c_i \delta_{\{x=1\}} = 0 \text{ for } x \in (\bar{x}, 1], \quad c_i = \int_{[0,1]} p_{i,W}(t,x) dx, \quad (44)$$

$$\frac{\partial p_{i,W}}{\partial t} + \frac{\partial F_{i,W}}{\partial x} + \mu p_{i,W} - \lambda p_{i,N} + \left( \sum_{j \in M, j \neq i} v_{i,j} \right) p_{i,W} - \sum_{j \in M, j \neq i} v_{j,i} p_{j,W} = 0, \quad x \in [0, \bar{x}] \quad (45)$$

with

$$p_{i,W}(t,x) = 0 \text{ for } x \in (\bar{x}, 1] \text{ and } p_{i,N}(t,1) = 0. \quad (46)$$

Again, we expect that the FPE is numerically recovered with a sufficiently large time $t \geq 0$. The initial condition must be non-negative as well as conservative ($\sum_{i \in M} \int_D (p_{i,W} + p_{i,N}) dx = 1$ at $t = 0$). The extended FPE is discretized using a classical upwind scheme combined with the WENO reconstruction [44]. A careful spatial discretization is required to maintain the mass conservation property. In addition, numerically solving FPEs with singular solutions like ours requires a careful discretization to correctly reproduce the solutions [53]. The discretization procedure of the numerical scheme for the system (43)-(46) is explained in **Appendix C**, showing that our scheme is conservative. The numerical experiment below demonstrates that the scheme correctly reproduces the singular solution.

### 4.1.3 Performance of the schemes

The two numerical schemes are examined with the exact solutions derived for the single-regime case before their applications to a real-world problem. The domain $D$ is discretized with $50n+1$ ($n = 1, 2, 4, 8, 16$) vertices and the convergence of numerical solutions to the exact solutions are evaluated. The parameter values are set as $\lambda = 1/7$, $\mu = 1$, $S = 0.07$, $c = 0.30$, $d = 0.20$, and $\delta = 0.1$. The threshold value $\bar{x}$ calculated from the specified parameter values is 0.807182, which was calculated using a trial and error method with the error smaller than $10^{-10}$. This value is utilized in the numerical computation of the FPE. Since the computation here focuses on the single-regime case, the terms on the regime-switching (terms containing $v_{i,j}$) are dropped from the equations. This simplification is not problematic because the source of the errors would be the first-order partial differential terms with respect to $x$, the coupled nature of the HJBE, and the singularity at $x=0$.

The time increment $\Delta t$ is set as $\Delta t = 5(\Delta x)^{1.5}$ for both the HJB and Fokker–Planck equations. This time increment is theoretically optimal for solving the HJBE, and we should have the order of 1.5 for the numerical solutions to the equation [43]. This means that computational errors measured in an error norm at the different resolutions $n=n_1$ and $n=n_2$ ($n_2 > n_1$), denoted as $e_1$ and $e_2$ respectively, should satisfy a scaling law not worse than $e_2 = 2^{-C(n_2-n_1)} e_1$ with the convergence rate $C=1.5$. For the FPE, it seems to be difficult to obtain a similar estimate because of the singularity. In the computation below, the threshold $\bar{x}$ is detected at a mid-point of some successive vertices. Each computation is terminated at the time step where the absolute difference between the numerical solutions at the successive time steps becomes smaller than the small error tolerance $10^{-14}$.

**Tables 1**-**2** show the $l^\infty$ errors (maximum difference between the numerical and exact solutions among all the vertices) between the numerical and exact $\Phi$ and the errors between the numerical and exact $\bar{x}$, respectively. The estimated convergence rate ($C$ in the previous paragraph) are also presented in the tables, which are better than the theoretical one discussed above. The error of $\bar{x}$ is smaller than the distance between each successive vertex.

**Table 3** shows the $l^\infty$ errors between the numerical and exact $p_N, p_W$ of the FPE for $x>0$. For each computation, the error is calculated for both $p_N$ and $p_W$, and the larger one is chosen. The computational results show that the numerical solutions, although with some irregularity in the convergence rate, converge to the exact solution as the resolution becomes fine for $x>0$. The error at the vertex placed at $x=0$ is considered in a separate table. The exact values of the weights of the Dirac Deltas at $x=0$ are 0.12534 and 0.02089 for the non-waiting case and waiting case, respectively. **Table 4** presents the absolute error on the singular part at $x=0$, where the singular parts of the numerical solutions are estimated as $p_N/h$ and $p_W/h$ with $h$ the cell size (**Appendix C**). The error decreases as the resolution increases. Consequently, we demonstrated that the employed schemes work successfully.

**Table 1.** Computed errors measured by the $l^\infty$ error and the corresponding convergence rates: the HJBE.

| Resolution $n$ | 1 | 2 | 4 | 8 | 16 |
|---|---|---|---|---|---|
| Error | 1.37.E-02 | 3.43.E-03 | 8.70.E-04 | 2.30.E-04 | 7.00.E-05 |
| Convergence rate | 2.03 | 1.99 | 1.93 | 1.72 | |

**Table 2.** Computed threshold value $\bar{x}$ and the corresponding error.

| Resolution $n$ | 1 | 2 | 4 | 8 | 16 |
|---|---|---|---|---|---|
| Computed $\bar{x}$ | 0.81 | 0.805 | 0.8075 | 0.80625 | 0.806875 |
| Error | 2.82.E-03 | 2.18.E-03 | 3.18.E-04 | 9.32.E-04 | 3.07.E-04 |

**Table 3.** Computed errors measured by the $l^\infty$ error and the corresponding convergence rates: the FPE. The error at the vertex placed at $x=0$ is not considered here.

| Resolution $n$ | 1 | 2 | 4 | 8 | 16 |
|---|---|---|---|---|---|
| Error | 1.02.E-02 | 2.01.E-03 | 1.96.E-03 | 8.71.E-04 | 2.89.E-04 |
| Convergence rate | 2.38 | 0.04 | 1.17 | 1.59 | |

**Table 4.** Computed errors of the weight of the singular part at $x=0$

| Resolution $n$ | 1 | 2 | 4 | 8 | 16 |
|---|---|---|---|---|---|
| Non-waiting case | 2.44.E-03 | 1.51.E-03 | 6.20.E-04 | 4.05.E-04 | 1.88.E-04 |
| Waiting case | 4.07.E-04 | 2.52.E-04 | 1.03.E-04 | 6.76.E-05 | 3.13.E-05 |

## 4.2 Computational conditions: sediment storage management

### 4.2.1 Problem background

The problem considered here is a sediment (earth and soils) replenishment problem in a dam-downstream river [23, 46]. A common issue that dam-downstream river environments worldwide encounter is the sediment trapping by dams [54], with which the sediment transport from the upstream is critically reduced or even stopped. The reduced sediment transport triggers an unnatural river environmental condition where nuisance algae grow thickly due to the reduced physical disturbance [23, 55, 56]. As a mitigation policy against the sediment trapping, replenishment of earth and soils from outside the river has been carried out in several case studies [57, 58]. It has been experimentally found that replenishing sediment is indeed effective for restoring the degraded dam-downstream river environment. Optimization of the sediment replenishment is a problem that has been paid less attention from a theoretical side. Especially, the existing models including the previous models by the authors [23, 46] do not consider the execution delay. This is the motivation of applying the present model to this sediment replenishment problem.

The case considered here is a problem of O Dam, H River, Japan [23, 46] to avoid sediment depletion in the dam-downstream reach. In this river, the O Dam has been working from 2011 for multiple purposes including water resources supply and flood mitigation. The sediment trapping is currently an environmental issue like many other rivers having dams. The stakeholders including a local fishery, the government of Japan, and the local government have been discussing about this issue in recent years. In April 2020, they initiated a sediment replenishment project and experimentally replenished the sediment in this month, with the amount of 100 (m$^3$). The planning of sediment replenishment is still not fixed, but the executions/observation can be possible at most weekly. The execution delay from a decision-making to the replenishment would require at least one to a few days because the earth and soils serving as the sediment must be transported from outside the river. In this section, we consider the optimal sediment replenishment policy with different values of the delay intensity $\mu$.

### 4.2.2 Coefficients and parameters

The coefficients and parameters are specified to carry out the numerical computation. Assume that a sediment lump with the maximum volume normalized to be 1 is placed in a dam-downstream river reach. The coefficient $S$ is specified as $S(i,x) = S_i \chi_{\{x>0\}}$ with non-negative constants $\{S_i\}_{i \in M}$. Assume that the river flow discharge, which is denoted as $\{q_i\}_{i \in M}$, is specified for each flow regime. We assume $q_i > 0$ ($i \in M$) and consider $S_i$ as a function of $q_i$ as $S(q_i) = S(q)|_{q=q_i}$. In this setting, $S$ physically means the transport rate of the sediment lump per unit time. A key fact is that $S$ can be estimated using the semi-empirical hydraulic formulae (Yoshioka [46]; Chapter 1 of Szymkiewicz [59]):

$$S(q) = \frac{1}{Y} \times 8 B \gamma^{1.5} \sqrt{g\sigma} \max\{\Theta - \Theta_c, 0\}^{\frac{3}{2}}, \quad \sigma = \frac{\rho_s}{\rho} - 1 > 0, \quad \Theta = \frac{\tau(q)}{\rho \sigma g \gamma}, \text{ and } \Theta_c = 0.047, \quad (47)$$

which is by the Meyer–Peter–Müller formula [60, 61] with

$$\tau(q) = \rho g h l = \rho g n_{\mathrm{M}}^{\frac{3}{5}} l^{\frac{7}{10}} B^{-\frac{3}{5}} q^{\frac{3}{5}}, \qquad (48)$$

which is by the Manning formula. Here, $\bar{Y}$ is the maximum volume of the storable sediment in a physical space, $B$ is the river channel width, $g$ is the gravitational acceleration, $\gamma$ is the diameter of sediment particles, $\rho$ is the density of water, $\rho_s$ is the soil density, $\tau$ is the bed shear stress as a positive and increasing function of $q > 0$, $n_{\mathrm{M}}$ is the roughness coefficient and $l$ is the channel slope. The relationship (47) implies that the transport rate is determined from the physical quantities, and that the sediment transport occur only when the discharge is large enough such that $\Theta > \Theta_c$. The parameters for the system dynamics are set as follows: $g = 9.81$ (m/s$^2$), $B = 25$ (m), $l = 0.001$ (-), $n_{\mathrm{M}} = 0.03$ (m$^{1/3}$/s), $\rho = 1,000$ (kg/m$^3$), $\rho_s = 2,600$ (kg/m$^3$), $\gamma = 5.0 \times 10^{-3}$ (m), and $\bar{Y} = 100$ (m$^3$).

The Markov chain $\alpha$ identified from an hourly discharge data [46] is utilized in this paper, but can be statistically estimated in any rivers if time-series of river flow discharge is available. The identification method itself is not the interest of this paper, but can be found in Yoshioka [46]. Yoshioka et al. [23] addressed a similar but coarser identification problem of a Markov chain from a discharge time-series data. We note that modeling river flow regimes by Markov chains is a common approach in hydrology and related research areas [62]. We set $I = 42$ with the discharge for each flow regime as an increasing function with respect to $i$: $q_i = 1.25 + 2.5i$ [46]. **Figure 4** shows the switching rates matrix $\nu = \left[ \nu_{i,j} \right]_{0 \le i,j \le 42}$ (1/day). Under this setting, we have $S_i = 0$ for $i = 0, 1$, meaning that the sediment transport does not occur during these low flow regimes. All the regimes are transient, and therefore the employed Markov chain satisfies the assumption (4). The stationary probability densities of the regimes $i = 0, 1, 2, 3, 4$ are 0.663, 0.124, 0.0962, 0.0378, 0.0140, respectively. The stationary probability densities for the other regimes are smaller than 0.01.

The other parameters are specified as follows. The discount rate $\delta$, which is the inverse of the time-scale of the decision-making, is set as 0.2 (1/day). By its definition, this corresponds to a decision-maker who has a sub-weekly perspective. Theoretically, we will be able to specify smaller value of $\delta$. We examined some of them but resulted in very slow (but stable) convergence that may not be practical. The parameters $c$ or $d$ are set as $c = 0.1$ and $d = 0.05$. Using a too large $c$ or $d$ seem to give a trivial policy that does not supply sediment at all (always supply the sediment), while too small values result in another trivial policy unconditionally replenishing the sediment at each observation. We set $\lambda = 1/7$ (1/day) assuming a weekly observation on average. Unless otherwise specified, we set $\mu = 1$ (1/day) assuming one-day delay on average.

On the computational resolution, we discretize $D$ with 350 and 175 equidistant vertices for the HJBE and FPE, respectively. The time increment for the HJBE is set as $5(1/350)^{1.5} = 7.64 \times 10^{-4}$ (day) and that for the FPE as $0.01/175 = 5.71 \times 10^{-5}$ (day). These values have been employed because using larger time increments was found to affect stability and accuracy of the numerical solutions. Each computation is terminated under the same criterion with that of Section 4.1.3.

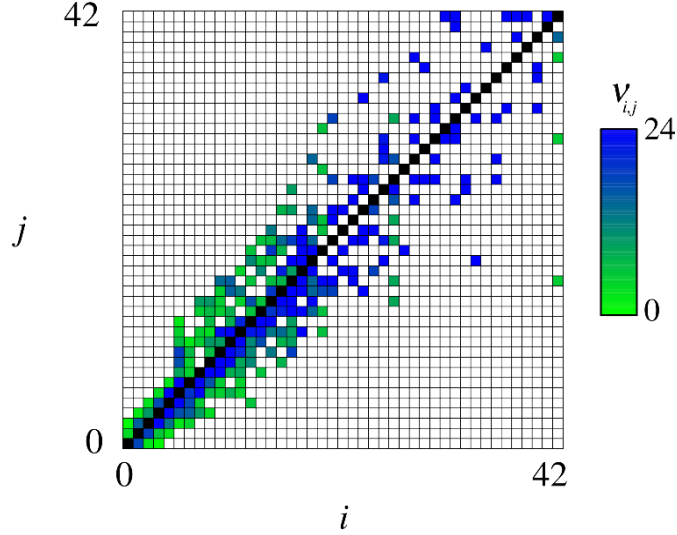

**Figure 4.** The switching rates matrix $\nu = \left[ \nu_{i,j} \right]_{0 \leq i, j \leq 42}$ (1/day) of the Markov chain. Notice that the diagonal elements (Colored black) are not used in our computation.

### 4.3 Impacts of the execution delay

We firstly present the computational results with $\lambda = 1/7$ (1/day) and $\mu = 1$ (1/day). **Figures 5-6** are representative figures of the computed optimal control $\eta = \eta^*(i, x)$ and value function $\Phi = \Phi(i, x)$ (**Figure 5**), and the optimal control $\eta = \eta^*(i, x)$ and the PDF $p = p(i, x)$ (**Figure 6**).

From **Figure 5**, the optimal control is almost activated ($\eta^*(i, x) = 1 - x$) for relatively small $i$. The non-monotone bang-bang nature of the control is due to utilizing the identified, and thus a non-parametric Markov chain estimated from an available data. A more regular profile of the optimal control may be obtained if we employ a parameterized model for the river flow regimes. Such models would be more tractable, but is possibly less flexible. The computational results suggest that the threshold type control considered in the tractable case applies to the flow regimes with relatively small $i$. The obtained numerical solution suggests that the sediment replenishment should not be carried out when the flow discharge is high, because more rapid sediment transport would occur in these regimes. Although not explicitly considered in this paper, replenishing the sediment under such high flow discharges may not be feasible in the real world. Considering this fact in the model would require the implementations cost to be regime-dependent and to specify a higher cost for larger $i$.

**Figure 6** demonstrates that the PDF is large for the small $i$ or along $x = 0$. The boundary singularity emerging in the exact solution in Section 3 also appears in the present case, implying that the solution is useful for analyzing the boundary behavior of the controlled dynamics. In **Figure 6**, we also find that the model has another singularity along the other boundary $x=1$ for the regimes $i=0,1$. Under the employed computational condition, the transport rate $S$ vanishes for these two regimes because of the small flow discharge such that the sediment particles are not transported ($\Theta \leq \Theta_c$). In addition, the identified Markov chain has a stationary PDF concentrated at the low flow regimes.

Based on these findings, the singular behavior of the PDF can be explained as follows. Once the sediment is replenished, the sediment storage becomes the full $x=1$. The sediment replenishment is carried out at the low flow regimes ($i=0,1$) under which with a high probability. Because the sediment transport does not occur in these regimes, the maximum value 1 persists for some finite interval until the flow regime transits to a relatively higher one ($i \geq 2$). This mechanism is not considered in the exact solution derived in Section 3 because it considers only one regime. Therefore, the singularity found along the boundary $x=1$ is a unique phenomenon emerging in this truly regime-switching model. The probabilities of the state along the boundaries $x=0$ and $x=1$ are estimated from the PDF as 0.209 and 0.231, respectively. Therefore, the sediment storage is either empty or full with the probability of 0.440, implying that evaluating the boundary states are important for understanding the controlled sediment storage dynamics subject to the delayed execution.

Finally, the impacts of execution delay are investigated from the two viewpoints suggested by the analysis above: the optimal control and the boundary behavior of the PDF. **Figure 7** compares the optimal controls for different values of $\mu$ (1/day), showing that the areas with positive $\eta^*$ increase as $\mu$ decreases. Recall that a smaller $\mu$ is likely to induce a larger execution delay. Due to the exponential distributional nature of the execution delay, the result implies that the optimal policy becomes getting closer to the policy with the unconditional sediment replenishment ($\eta^* > 0$ except along $x=1$). **Figure 8** compares the boundary probabilities $P(0) = \Pr(x=0)$, $P(1) = \Pr(x=1)$, and $P(0) + P(1)$ for different values of the delay intensity $\mu$ (1/day) under the corresponding optimal controls. The computational results show that the boundary probabilities are in general decreasing with respect to $\mu$. Especially, the probability $P(0) + P(1)$ rapidly increases as $\mu$ becomes smaller than 1. Therefore, there would exist a critical value of $\mu$, below which (resp., above which) the boundary probability is sensitive (resp., insensitive to) the execution delay. In our case, the execution delay larger than one day may on average significantly increase the probabilities that the sediment storage is either empty or full at the observations. In addition, the computational results suggest that the controlled dynamics would not be critically changed by forcing the delay intensity $\mu$ to become larger than 1 (1/day).

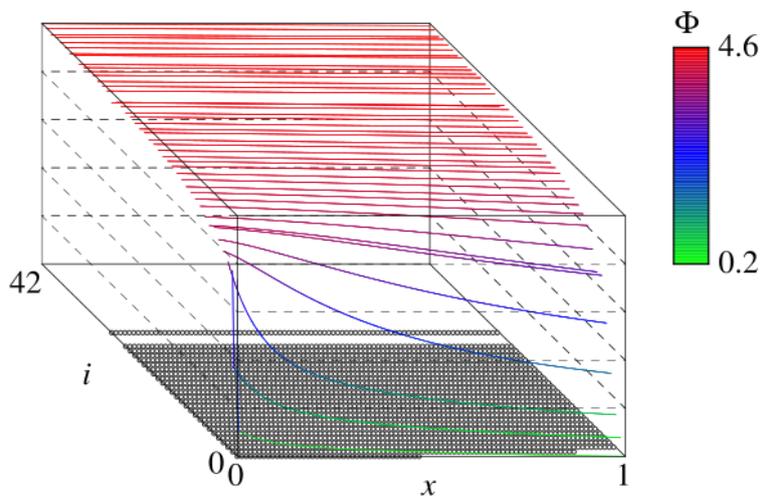

**Figure 5.** The computed value function $\Phi$ (Colored curves) and the associated optimal control $\eta^*$ (Positive where the circles present).

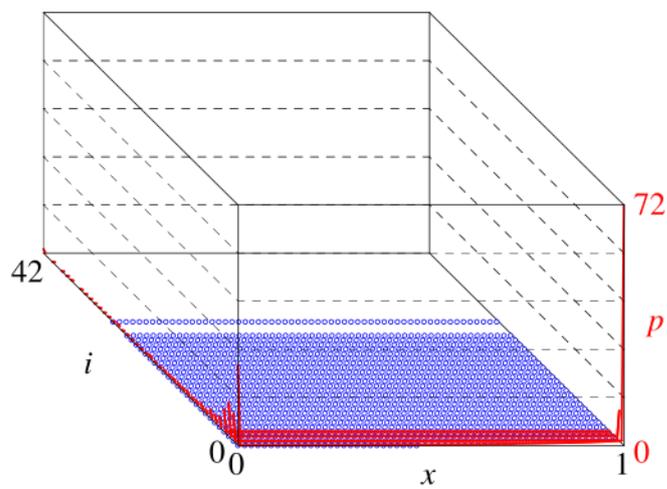

**Figure 6.** The computed PDF $p$ (Red curves) and the associated optimal control $\eta^*$ (Positive where the circles present). The PDF is plotted only at the points where $p \geq 10^{-4}$.

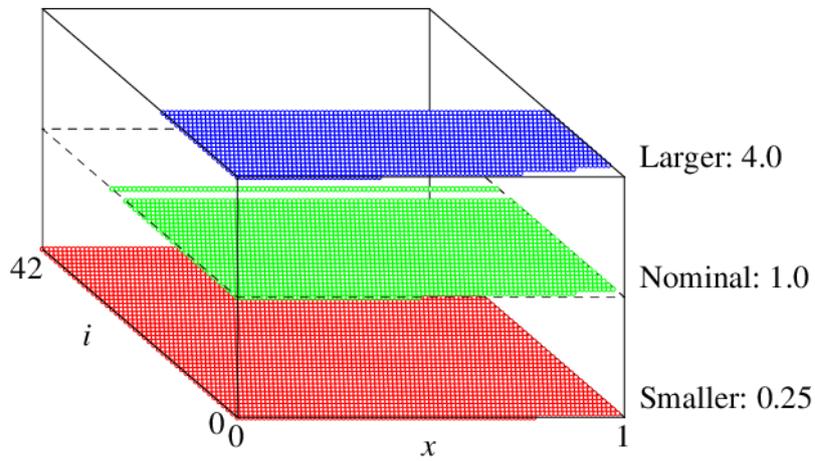

**Figure 7.** Comparison of the optimal controls for different values of the delay intensity $\mu$ (1/day) (Blue: larger value $\mu = 4.0$, Green: nominal value $\mu = 1.0$, and Red: smaller value $\mu = 0.25$). The areas with $\eta^* > 0$ are plotted in the figure. For $\mu = 0.25$, the computed $\eta^*$ is positive for positive $i$

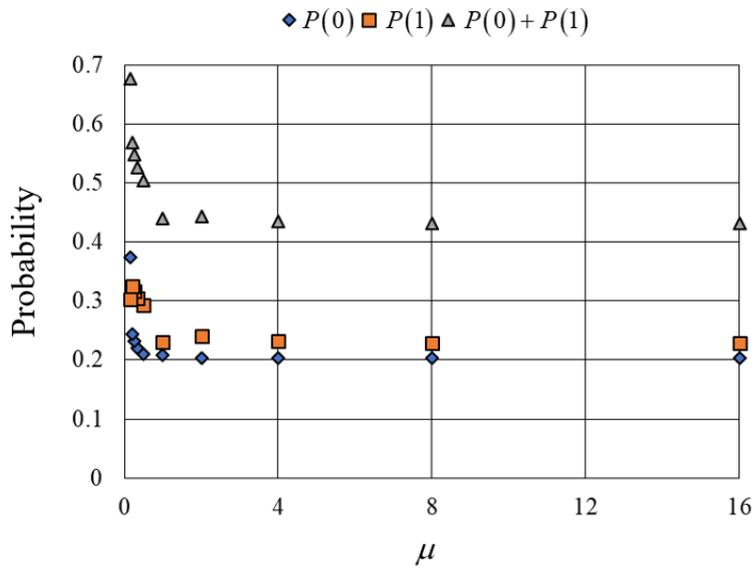

**Figure 8.** Comparison of the boundary probabilities $P(0)$, $P(1)$, and $P(0) + P(1)$ for different values of the delay intensity $\mu$ (1/day).

## 5. Conclusions

We considered a stochastic control problem of a non-smooth dynamics with discrete random observations and delayed impulsive interventions. We presented the HJBE governing the optimal control and the value function, and the FPE governing the PDF of the controlled dynamics. We demonstrated that they are tractable under a simplified case. The HJBE under the simplification is verified to admit a threshold type optimal policy. The FPE under the simplification admits a PDF having a boundary singularity. The latter was validated with a Monte-Carlo method. We utilized numerical schemes for discretizing these equations and demonstrated their convergence. The two equations were finally applied to a recent sediment storage dynamics problem in a river, implying an importance of considering the execution delay in the application.

We could numerically analyze the delay-dependence of the optimal control and the controlled PDF; however, their mathematical justification is still an open issue. The problem considered in this paper assumed an execution delay following an exponential distribution; however, more general distributions, such as a gamma distribution, may be more realistic. In such a case, both the HJBE and the FPE would become more complicated. The computational results in Section 4 suggest a possibility to more efficiently model the rapid flushing out of the sediment using a jump process. We will approach this issue by randomizing the term $S\mathrm{d}t$ using Lévy processes. Problems with a seasonality can be handled with the presented framework by considering an appropriate finite-horizon problem, but would require a more careful treatment of the value function and the HJBE especially near the terminal time [34, 35]. A viscosity solution approach would be necessary for handling this issue in both theory [63] and computation [64]. Analyzing impacts of delayed execution in the other non-smooth dynamics arising from engineering problems, such as fishery resources dynamics [65] and ecological dynamics [66], is an important research topic from both theoretical and practical viewpoints. Higher-dimensional problems are also interesting topics, which may be approached by utilizing backward SDEs [67]; however, problems with the discrete and random observations with delayed impulsive interventions seem to be least explored. Development of a fast and accurate numerical scheme that is able to efficiently handle small $\delta$ is also a key challenge for both the proposed and extended problems.

## Appendices
**Appendix A: Technical proofs**

Proofs of **Propositions 1** and **2** are presented.

**(Proof of Proposition 1)**

Recall that $\eta_0 = 0$ and that $\tau_k$ is the $k$ th observation time at which the decision is made, while the execution is carried out at $\tau_k + \omega_k$ if $\eta_k > 0$. Choose some $x \in D$ and $\bar{\eta} \in C$. Hereafter, $\mathbb{E}^{X_s}$ represents the expectation conditioned on the time $s$ and the value $X_s$, and $\mathbb{E}^{t,X_s}$ that conditioned on time $t$ and the value $X_s$. Therefore, $\mathbb{E}^{s,X_s} = \mathbb{E}^{X_s}$ and $\mathbb{E}^{X_0} = \mathbb{E}^x$. Recall that the controlled dynamics are autonomous.

By $\bar{\Phi} \in C^1(D)$, for each $k \geq 0$, applying Itô's formula to $e^{-\delta t}\bar{\Phi}(X_t)$ yields

$$e^{-\delta(\tau_k+\Delta_k)}\bar{\Phi}(X_{\tau_k+\Delta_k}) - e^{-\delta\tau_k}\bar{\Phi}(X_{\tau_k}) = \int_{\tau_k}^{\tau_k+\Delta_k} e^{-\delta s}\left(-\delta\bar{\Phi}(X_s) - S\chi_{\{X_s>0\}}\frac{d\bar{\Phi}}{dx}(X_s)\right)ds \quad (49)$$

and

$$e^{-\delta\tau_{k+1}}\bar{\Phi}(X_{\tau_{k+1}}) - e^{-\delta(\tau_k+\Delta_k)}\bar{\Phi}(X_{(\tau_k+\Delta_k)+}) = \int_{\tau_k+\Delta_k}^{\tau_{k+1}} e^{-\delta s}\left(-\delta\bar{\Phi}(X_s) - S\chi_{\{X_s>0\}}\frac{d\bar{\Phi}}{dx}(X_s)\right)ds. \quad (50)$$

Combining (49) and (50) yields

$$\begin{aligned}e^{-\delta\tau_{k+1}}\bar{\Phi}(X_{\tau_{k+1}}) - e^{-\delta\tau_k}\bar{\Phi}(X_{\tau_k}) &= \int_{\tau_k}^{\tau_{k+1}} e^{-\delta s}\left(-\delta\bar{\Phi}(X_s) - S\chi_{\{X_s>0\}}\frac{d\bar{\Phi}}{dx}(X_s)\right)ds \\ &+ e^{-\delta(\tau_k+\Delta_k)}\left[\bar{\Phi}(X_{(\tau_k+\Delta_k)+}) - \bar{\Phi}(X_{\tau_k+\Delta_k})\right]\end{aligned} \quad (51)$$

Since $\bar{\Phi}$ solves (17) pointwise, we have

$$\begin{aligned}e^{-\delta\tau_{k+1}}\bar{\Phi}(X_{\tau_{k+1}}) - e^{-\delta\tau_k}\bar{\Phi}(X_{\tau_k}) &= \int_{\tau_k}^{\tau_{k+1}} e^{-\delta s}\left(\lambda\left(\bar{\Phi}(X_s) - \min\{\bar{\Phi}(X_s), \hat{\bar{\Phi}}(X_s)\}\right) - \chi_{\{X_s=0\}}\right)ds \\ &+ e^{-\delta(\tau_k+\Delta_k)}\left[\bar{\Phi}(X_{(\tau_k+\Delta_k)+}) - \bar{\Phi}(X_{\tau_k+\Delta_k})\right]\end{aligned} \quad (52)$$

We then get

$$\begin{aligned}&\int_{\tau_k}^{\tau_{k+1}} e^{-\delta s}\chi_{\{X_s=0\}} ds + e^{-\delta(\tau_k+\Delta_k)}\left(c\eta_k + d\chi_{\{\eta_k>0\}}\right) + e^{-\delta\tau_{k+1}}\bar{\Phi}(X_{\tau_{k+1}}) - e^{-\delta\tau_k}\bar{\Phi}(X_{\tau_k}) \\ &= \int_{\tau_k}^{\tau_{k+1}} e^{-\delta s}\lambda\left(\bar{\Phi}(X_s) - \min\{\bar{\Phi}(X_s), \hat{\bar{\Phi}}(X_s)\}\right)ds \\ &+ e^{-\delta(\tau_k+\Delta_k)}\left(\bar{\Phi}(X_{(\tau_k+\Delta_k)+}) - \bar{\Phi}(X_{\tau_k+\Delta_k}) + c\eta_k + d\chi_{\{\eta_k>0\}}\right)\end{aligned} \quad (53)$$

Taking the expectation $\mathbb{E}^x$ and summing up for $k = 0, 1, 2, \ldots$ yields (Notice that $\tau_k \to +\infty$ a.s. for $k \to +\infty$, $\tau_0 = 0$, $\eta_0 = 0$, $\bar{\Phi} \in C^1(D)$, and $D$ is compact. We can exchange the expectation and summations in the equations below since each term is bounded.)

$$\begin{aligned}&\sum_{k=0}^{\infty}\mathbb{E}^x\left[\int_{\tau_k}^{\tau_{k+1}} e^{-\delta s}\chi_{\{X_s=0\}} ds + e^{-\delta(\tau_k+\Delta_k)}\left(c\eta_k + d\chi_{\{\eta_k>0\}}\right)\right] \\ &= \mathbb{E}^x\left[\int_0^{\infty} e^{-\delta s}\chi_{\{X_s=0\}} ds + \sum_{k=1}^{\infty} e^{-\delta(\tau_k+\Delta_k)}\left(c\eta_k + d\chi_{\{\eta_k>0\}}\right)\right], \\ &= \phi(x, \bar{\eta})\end{aligned} \quad (54)$$

$$\sum_{k=0}^{\infty}\mathbb{E}^x\left[e^{-\delta\tau_{k+1}}\bar{\Phi}(X_{\tau_{k+1}}) - e^{-\delta\tau_k}\bar{\Phi}(X_{\tau_k})\right] = -\bar{\Phi}(X_{\tau_0}) = -\bar{\Phi}(x), \quad (55)$$

and thus

$$\begin{aligned}&\phi(x,\bar{\eta}) - \bar{\Phi}(x) \\ &= \sum_{k=0}^{\infty}\mathbb{E}^x\left[\int_{\tau_k}^{\tau_{k+1}} e^{-\delta s}\lambda\left(\bar{\Phi}(X_s) - \min\{\bar{\Phi}(X_s), \hat{\bar{\Phi}}(X_s)\}\right)ds\right] \\ &+ \sum_{k=0}^{\infty}\mathbb{E}^x\left[e^{-\delta(\tau_k+\Delta_k)}\left(\bar{\Phi}(X_{(\tau_k+\Delta_k)+}) - \bar{\Phi}(X_{\tau_k+\Delta_k}) + c\eta_k + d\chi_{\{\eta_k>0\}}\right)\right]\end{aligned} \quad (56)$$

Now, we have

$$\mathbb{E}^x\left[\chi_{\{\eta_k=0\}}e^{-\delta(\tau_k+\Delta_k)}\left(\bar{\Phi}(X_{(\tau_k+\Delta_k)+})-\bar{\Phi}(X_{\tau_k+\Delta_k})+c\eta_k+d\chi_{\{\eta_k>0\}}\right)\right]=0 \tag{57}$$

by $X_{(\tau_k+\Delta_k)+}=X_{\tau_k+\Delta_k}$ if $\eta_k=0$. On the other hand, if $\eta_k\left(=1-X_{\tau_k+\Delta_k}\right)>0$, then $\Delta_k=\omega_k$, $X_{(\tau_k+\Delta_k)+}=1$, $\eta_k=1-X_{\tau_k+\omega_k}$, and by adding and subtracting an integral term, we get

$$\begin{aligned}
&\mathbb{E}^x\left[\chi_{\{\eta_k>0\}}e^{-\delta(\tau_k+\Delta_k)}\left(\bar{\Phi}(X_{(\tau_k+\Delta_k)+})-\bar{\Phi}(X_{\tau_k+\Delta_k})+c\eta_k+d\chi_{\{\eta_k>0\}}\right)\right] \\
&=\mathbb{E}^x\left[\chi_{\{\eta_k>0\}}e^{-\delta\tau_k}\mathbb{E}^{X_{\tau_k}}\left[e^{-\delta\omega_k}\left(\bar{\Phi}(1)+c(1-X_{\tau_k+\omega_k})+d\right)+\int_{\tau_k}^{\tau_k+\omega_k}e^{-\delta(s-\tau_k)}\chi_{\{X_s=0\}}\mathrm{d}s\right]\right] \\
&\quad-\mathbb{E}^x\left[\chi_{\{\eta_k>0\}}\left(e^{-\delta(\tau_k+\Delta_k)}\bar{\Phi}(X_{\tau_k+\Delta_k})+\int_{\tau_k}^{\tau_k+\Delta_k}e^{-\delta s}\chi_{\{X_s=0\}}\mathrm{d}s\right)\right] \\
&=\mathbb{E}^x\left[\chi_{\{\eta_k>0\}}e^{-\delta\tau_k}\mathbb{E}^{0,X_{\tau_k}}\left[e^{-\delta\omega_k}\left(\bar{\Phi}(1)+c(1-X_{\omega_k})+d\right)+\int_0^{\omega_k}e^{-\delta s}\chi_{\{X_s=0\}}\mathrm{d}s\right]\right] \\
&\quad-\mathbb{E}^x\left[\chi_{\{\eta_k>0\}}\left(e^{-\delta(\tau_k+\Delta_k)}\bar{\Phi}(X_{\tau_k+\Delta_k})+\int_{\tau_k}^{\tau_k+\Delta_k}e^{-\delta s}\chi_{\{X_s=0\}}\mathrm{d}s\right)\right] \\
&=\mathbb{E}^x\left[\chi_{\{\eta_k>0\}}e^{-\delta\tau_k}\hat{\bar{\Phi}}(X_{\tau_k})\right]-\mathbb{E}^x\left[\chi_{\{\eta_k>0\}}\left(e^{-\delta(\tau_k+\Delta_k)}\bar{\Phi}(X_{\tau_k+\Delta_k})+\int_{\tau_k}^{\tau_k+\Delta_k}e^{-\delta s}\chi_{\{X_s=0\}}\mathrm{d}s\right)\right]
\end{aligned} \tag{58}$$

By (57) and (58), we get

$$\begin{aligned}
&\mathbb{E}^x\left[e^{-\delta(\tau_k+\Delta_k)}\left(\bar{\Phi}(X_{(\tau_k+\Delta_k)+})-\bar{\Phi}(X_{\tau_k+\Delta_k})+c\eta_k+d\chi_{\{\eta_k>0\}}\right)\right] \\
&=\mathbb{E}^x\left[\chi_{\{\eta_k>0\}}e^{-\delta\tau_k}\hat{\bar{\Phi}}(X_{\tau_k})\right]-\mathbb{E}^x\left[\chi_{\{\eta_k>0\}}\left(e^{-\delta(\tau_k+\Delta_k)}\bar{\Phi}(X_{\tau_k+\Delta_k})+\int_{\tau_k}^{\tau_k+\Delta_k}e^{-\delta s}\chi_{\{X_s=0\}}\mathrm{d}s\right)\right]
\end{aligned} \tag{59}$$

By (56) and (59), we arrive at

$$\begin{aligned}
&\phi(x,\bar{\eta})-\bar{\Phi}(x) \\
&=\sum_{k=0}\mathbb{E}^x\left[\int_{\tau_k}^{\tau_{k+1}}e^{-\delta s}\lambda\left(\bar{\Phi}(X_s)-\min\{\bar{\Phi}(X_s),\hat{\bar{\Phi}}(X_s)\}\right)\mathrm{d}s\right] \\
&\quad+\sum_{k=0}\left(\mathbb{E}^x\left[\chi_{\{\eta_k>0\}}e^{-\delta\tau_k}\hat{\bar{\Phi}}(X_{\tau_k})\right]-\mathbb{E}^x\left[\chi_{\{\eta_k>0\}}e^{-\delta(\tau_k+\Delta_k)}\bar{\Phi}(X_{\tau_k+\Delta_k})+\int_{\tau_k}^{\tau_k+\Delta_k}e^{-\delta s}\chi_{\{X_s=0\}}\mathrm{d}s\right]\right)
\end{aligned} \tag{60}$$

Again using Itô's formula combined with (17) yields

$$\begin{aligned}
&e^{-\delta(\tau_k+\Delta_k)}\bar{\Phi}(X_{\tau_k+\Delta_k}) \\
&=e^{-\delta\tau_k}\bar{\Phi}(X_{\tau_k})+\int_{\tau_k}^{\tau_k+\Delta_k}e^{-\delta s}\lambda\left(\bar{\Phi}(X_s)-\min\{\bar{\Phi}(X_s),\hat{\bar{\Phi}}(X_s)\}\right)\mathrm{d}s-\int_{\tau_k}^{\tau_k+\Delta_k}e^{-\delta s}\chi_{\{X_s=0\}}\mathrm{d}s
\end{aligned} \tag{61}$$

Substituting (61) into (60) yields

$$\phi(x,\bar{\eta}) - \bar{\Phi}(x)$$
$$= \sum_{k=0}^{\infty} \mathbb{E}^x \left[ \int_{\tau_k}^{\tau_{k+1}} e^{-\delta s} \lambda \left( \bar{\Phi}(X_s) - \min\{\bar{\Phi}(X_s), \hat{\bar{\Phi}}(X_s)\} \right) ds \right]$$
$$+ \sum_{k=0}^{\infty} \left( \mathbb{E}^x \left[ \chi_{\{\eta_k > 0\}} e^{-\delta \tau_k} \hat{\bar{\Phi}}(X_{\tau_k}) \right] - \mathbb{E}^x \left[ \begin{array}{l} \chi_{\{\eta_k > 0\}} e^{-\delta \tau_k} \bar{\Phi}(X_{\tau_k}) \\ + \int_{\tau_k}^{\tau_k + \Delta_k} e^{-\delta s} \lambda \left( \bar{\Phi}(X_s) - \min\{\bar{\Phi}(X_s), \hat{\bar{\Phi}}(X_s)\} \right) ds \end{array} \right] \right)$$
$$= \sum_{k=0}^{\infty} \mathbb{E}^x \left[ \int_{\tau_k + \Delta_k}^{\tau_{k+1}} e^{-\delta s} \lambda \left( \bar{\Phi}(X_s) - \min\{\bar{\Phi}(X_s), \hat{\bar{\Phi}}(X_s)\} \right) ds \right] \quad , \quad (62)$$
$$+ \sum_{k=0}^{\infty} \mathbb{E}^x \left[ \chi_{\{\eta_k > 0\}} e^{-\delta \tau_k} \left( \hat{\bar{\Phi}}(X_{\tau_k}) - \bar{\Phi}(X_{\tau_k}) \right) \right]$$
$$= \sum_{k=0}^{\infty} \mathbb{E}^x \left[ \int_{\tau_k + \Delta_k}^{\tau_{k+1}} e^{-\delta s} \lambda \left( \bar{\Phi}(X_s) - \min\{\bar{\Phi}(X_s), \hat{\bar{\Phi}}(X_s)\} \right) ds \right]$$
$$+ \sum_{k=0}^{\infty} \mathbb{E}^x \left[ \chi_{\{\eta_{k+1} > 0\}} e^{-\delta \tau_{k+1}} \left( \hat{\bar{\Phi}}(X_{\tau_{k+1}}) - \bar{\Phi}(X_{\tau_{k+1}}) \right) \right]$$

where we used $\eta_0 = 0$ to derive the last line, and also the identity

$$\mathbb{E}^x \left[ \chi_{\{\eta_k > 0\}} \int_{\tau_k}^{\tau_k + \Delta_k} e^{-\delta s} \lambda \left( \bar{\Phi}(X_s) - \min\{\bar{\Phi}(X_s), \hat{\bar{\Phi}}(X_s)\} \right) ds \right]$$
$$= \mathbb{E}^x \left[ \chi_{\{\eta_k > 0\}} \int_{\tau_k}^{\tau_k + \Delta_k} e^{-\delta s} \lambda \left( \bar{\Phi}(X_s) - \min\{\bar{\Phi}(X_s), \hat{\bar{\Phi}}(X_s)\} \right) ds \right]$$
$$+ \mathbb{E}^x \left[ \chi_{\{\eta_k = 0\}} \int_{\tau_k}^{\tau_k + \Delta_k} e^{-\delta s} \lambda \left( \bar{\Phi}(X_s) - \min\{\bar{\Phi}(X_s), \hat{\bar{\Phi}}(X_s)\} \right) ds \right] \quad (63)$$
$$= \mathbb{E}^x \left[ \int_{\tau_k}^{\tau_k + \Delta_k} e^{-\delta s} \lambda \left( \bar{\Phi}(X_s) - \min\{\bar{\Phi}(X_s), \hat{\bar{\Phi}}(X_s)\} \right) ds \right]$$

because of $\Delta_k = 0$ if $\eta_k = 0$ (The third expectation equals 0). Now, we have

$$\mathbb{E}^x \left[ \int_{\tau_k + \Delta_k}^{\tau_{k+1}} e^{-\delta s} \lambda \left( \bar{\Phi}(X_s) - \min\{\bar{\Phi}(X_s), \hat{\bar{\Phi}}(X_s)\} \right) ds \right]$$
$$= \mathbb{E}^x \left[ \mathbb{E}^{X_{\tau_k + \Delta_k}} \left[ \int_{\tau_k + \Delta_k}^{\tau_{k+1}} e^{-\delta s} \lambda \left( \bar{\Phi}(X_s) - \min\{\bar{\Phi}(X_s), \hat{\bar{\Phi}}(X_s)\} \right) ds \right] \right]$$
$$= \mathbb{E}^x \left[ \mathbb{E}^{X_{\tau_k + \Delta_k}} \left[ \int_0^{\tau_{k+1} - (\tau_k + \Delta_k)} e^{-\delta(s + (\tau_k + \Delta_k))} \lambda \left( \bar{\Phi}(X_{s+(\tau_k + \Delta_k)}) - \min\{\bar{\Phi}(X_{s+(\tau_k + \Delta_k)}), \hat{\bar{\Phi}}(X_{s+(\tau_k + \Delta_k)})\} \right) ds \right] \right]$$
$$= \mathbb{E}^x \left[ \mathbb{E}^{0, X_{\tau_k + \Delta_k}} \left[ \int_0^{\tau_{k+1} - (\tau_k + \Delta_k)} e^{-\delta(s + (\tau_k + \Delta_k))} \lambda \left( \bar{\Phi}(X_s) - \min\{\bar{\Phi}(X_s), \hat{\bar{\Phi}}(X_s)\} \right) ds \right] \right] \quad (64)$$
$$= \mathbb{E}^x \left[ \lambda e^{-\delta(\tau_k + \Delta_k)} \mathbb{E}^{0, X_{\tau_k + \Delta_k}} \left[ \int_0^{\tau_{k+1} - (\tau_k + \Delta_k)} e^{-\delta s} \left( \bar{\Phi}(X_s) - \min\{\bar{\Phi}(X_s), \hat{\bar{\Phi}}(X_s)\} \right) ds \right] \right]$$
$$= \mathbb{E}^x \left[ \lambda e^{-\delta(\tau_k + \Delta_k)} \mathbb{E}^{0, X_{\tau_k + \Delta_k}} \left[ \int_0^{+\infty} e^{-(\delta + \lambda)s} \left( \bar{\Phi}(X_s) - \min\{\bar{\Phi}(X_s), \hat{\bar{\Phi}}(X_s)\} \right) ds \right] \right]$$

by $\tau_{k+1} - (\tau_k + \Delta_k) \sim \text{Exp}_\lambda$. Since $\chi_{\{\eta_{k+1} > 0\}} = 0$ or 1, we have the inequality

$$\mathbb{E}^x\left[\chi_{\{\eta_{k+1}>0\}}e^{-\delta\tau_{k+1}}\left(\hat{\bar{\Phi}}(X_{\tau_{k+1}})-\bar{\Phi}(X_{\tau_{k+1}})\right)\right]$$

$$\geq \mathbb{E}^x\left[\chi_{\{\eta_{k+1}>0\}}e^{-\delta\tau_{k+1}}\left(\min\left\{\bar{\Phi}(X_{\tau_{k+1}}),\hat{\bar{\Phi}}(X_{\tau_{k+1}})\right\}-\bar{\Phi}(X_{\tau_{k+1}})\right)\right]$$

$$=-\mathbb{E}^x\left[\chi_{\{\eta_{k+1}>0\}}e^{-\delta\tau_{k+1}}\left(\bar{\Phi}(X_{\tau_{k+1}})-\min\left\{\bar{\Phi}(X_{\tau_{k+1}}),\hat{\bar{\Phi}}(X_{\tau_{k+1}})\right\}\right)\right] \quad (65)$$

$$\geq -\mathbb{E}^x\left[e^{-\delta\tau_{k+1}}\left(\bar{\Phi}(X_{\tau_{k+1}})-\min\left\{\bar{\Phi}(X_{\tau_{k+1}}),\hat{\bar{\Phi}}(X_{\tau_{k+1}})\right\}\right)\right]$$

by $\bar{\Phi}(X_{\tau_{k+1}})-\min\{\hat{\bar{\Phi}}(X_{\tau_{k+1}}),\hat{\bar{\Phi}}(X_{\tau_{k+1}})\}\geq 0$. This inequality can be further calculated as

$$-\mathbb{E}^x\left[e^{-\delta\tau_{k+1}}\left(\bar{\Phi}(X_{\tau_{k+1}})-\min\left\{\bar{\Phi}(X_{\tau_{k+1}}),\hat{\bar{\Phi}}(X_{\tau_{k+1}})\right\}\right)\right]$$

$$=-\mathbb{E}^x\left[\mathbb{E}^{X_{\tau_k+\Delta_k}}\left[e^{-\delta\tau_{k+1}}\left(\bar{\Phi}(X_{\tau_{k+1}})-\min\left\{\bar{\Phi}(X_{\tau_{k+1}})\hat{\bar{\Phi}}(X_{\tau_{k+1}})\right\}\right)\right]\right]$$

$$=-\mathbb{E}^x\left[e^{-\delta(\tau_k+\Delta_k)}\mathbb{E}^{X_{\tau_k+\Delta_k}}\left[e^{-\delta(\tau_{k+1}-(\tau_k+\Delta_k))}\left(\bar{\Phi}(X_{\tau_{k+1}})-\min\left\{\bar{\Phi}(X_{\tau_{k+1}})\hat{\bar{\Phi}}(X_{\tau_{k+1}})\right\}\right)\right]\right]$$

$$=-\mathbb{E}^x\left[e^{-\delta(\tau_k+\Delta_k)}\mathbb{E}^{0,X_{\tau_k+\Delta_k}}\left[e^{-\delta(\tau_{k+1}-(\tau_k+\Delta_k))}\left(\bar{\Phi}(X_{\tau_{k+1}-(\tau_k+\Delta_k)})-\min\left\{\bar{\Phi}(X_{\tau_{k+1}-(\tau_k+\Delta_k)}),\hat{\bar{\Phi}}(X_{\tau_{k+1}-(\tau_k+\Delta_k)})\right\}\right)\right]\right], \quad (66)$$

$$=-\mathbb{E}^x\left[e^{-\delta(\tau_k+\Delta_k)}\mathbb{E}^{0,X_{\tau_k+\Delta_k}}\left[\int_0^\infty \lambda e^{-\lambda s}e^{-\delta s}\left(\bar{\Phi}(X_s)-\min\left\{\bar{\Phi}(X_s),\hat{\bar{\Phi}}(X_s)\right\}\right)\mathrm{d}s\right]\right]$$

$$=-\mathbb{E}^x\left[\lambda e^{-\delta(\tau_k+\Delta_k)}\mathbb{E}^{0,X_{\tau_k+\Delta_k}}\left[\int_0^{+\infty} e^{-(\lambda+\delta)s}\left(\bar{\Phi}(X_s)-\min\left\{\bar{\Phi}(X_s),\hat{\bar{\Phi}}(X_s)\right\}\right)\mathrm{d}s\right]\right]$$

again by $\tau_{k+1}-(\tau_k+\Delta_k)\sim\mathrm{Exp}_\lambda$. By (64) and (66), (62) is evaluated from below as

$$\phi(x,\bar{\eta})-\bar{\Phi}(x)$$
$$=\sum_{k=0}\mathbb{E}^x\left[\int_{\tau_k+\Delta_k}^{\tau_{k+1}}e^{-\delta s}\lambda\left(\bar{\Phi}(X_s)-\min\left\{\bar{\Phi}(X_s),\hat{\bar{\Phi}}(X_s)\right\}\right)\mathrm{d}s\right]$$
$$+\sum_{k=0}\mathbb{E}^x\left[\chi_{\{\eta_{k+1}>0\}}e^{-\delta\tau_{k+1}}\left(\hat{\bar{\Phi}}(X_{\tau_{k+1}})-\bar{\Phi}(X_{\tau_{k+1}})\right)\right]$$
$$\geq\sum_{k=0}\mathbb{E}^x\left[\lambda e^{-\delta(\tau_k+\Delta_k)}\mathbb{E}^{0,X_{\tau_k+\Delta_k}}\left[\int_0^{+\infty}e^{-(\delta+\lambda)s}\left(\bar{\Phi}(X_s)-\min\left\{\bar{\Phi}(X_s),\hat{\bar{\Phi}}(X_s)\right\}\right)\mathrm{d}s\right]\right] \quad , \quad (67)$$
$$+\sum_{k=0}\left(-\mathbb{E}^x\left[\lambda e^{-\delta(\tau_k+\Delta_k)}\mathbb{E}^{0,X_{\tau_k+\Delta_k}}\left[\int_0^{+\infty}e^{-(\lambda+\delta)s}\left(\bar{\Phi}(X_s)-\min\left\{\bar{\Phi}(X_s),\hat{\bar{\Phi}}(X_s)\right\}\right)\mathrm{d}s\right]\right]\right)$$
$$=0$$

and thus $\phi(x,\bar{\eta})\geq\bar{\Phi}(x)$. Since $x\in D$ and $\bar{\eta}\in C$ are arbitrary, we get

$$\Phi\geq\bar{\Phi} \quad \text{in} \quad D. \quad (68)$$

The next step is to show $\Phi=\bar{\Phi}$. The control with (20) is denoted as $\bar{\eta}=\hat{\eta}$. This control, which is clearly admissible, is used in the rest of this proof. By (20), we have

$$\bar{\Phi}(X_s)-\min\left\{\bar{\Phi}(X_s),\hat{\bar{\Phi}}(X_s)\right\}=\begin{cases}0 & (X_s>\bar{x})\\ \bar{\Phi}(X_s)-\hat{\bar{\Phi}}(X_s) & (X_s\leq\bar{x})\end{cases}. \quad (69)$$

We again utilize (62) and (64). The last line of (64) is now evaluated as

$$\mathbb{E}^x\left[\lambda e^{-\delta(\tau_k+\Delta_k)}\mathbb{E}^{0,X_{\tau_k+\Delta_k}}\left[\int_0^{+\infty}e^{-(\delta+\lambda)s}\left(\bar{\Phi}(X_s)-\min\left\{\bar{\Phi}(X_s),\hat{\bar{\Phi}}(X_s)\right\}\right)\mathrm{d}s\right]\right]$$
$$=\mathbb{E}^x\left[\lambda e^{-\delta(\tau_k+\Delta_k)}\mathbb{E}^{0,X_{\tau_k+\Delta_k}}\left[\int_0^{+\infty}e^{-(\delta+\lambda)s}\chi_{\{X_s\leq\bar{x}\}}\left(\bar{\Phi}(X_s)-\hat{\bar{\Phi}}(X_s)\right)\mathrm{d}s\right]\right]$$ (70)

By (69), $\chi_{\{\eta_{k+1}\leq\bar{x}\}}=\chi_{\{X_{\tau_{k+1}}\leq\bar{x}\}}$ for $k\geq 0$. Therefore, a calculation like (65)-(66) yields (Notice that the result below is an identity but not an inequality)

$$\mathbb{E}^x\left[\chi_{\{\eta_{k+1}>0\}}e^{-\delta\tau_{k+1}}\left(\hat{\bar{\Phi}}(X_{\tau_{k+1}})-\bar{\Phi}(X_{\tau_{k+1}})\right)\right]$$
$$=\mathbb{E}^x\left[\chi_{\{X_{\tau_{k+1}}\leq\bar{x}\}}e^{-\delta\tau_{k+1}}\left(\hat{\bar{\Phi}}(X_{\tau_{k+1}})-\bar{\Phi}(X_{\tau_{k+1}})\right)\right]$$
$$=-\mathbb{E}^x\left[\chi_{\{X_{\tau_{k+1}}\leq\bar{x}\}}e^{-\delta\tau_{k+1}}\left(\bar{\Phi}(X_{\tau_{k+1}})-\hat{\bar{\Phi}}(X_{\tau_{k+1}})\right)\right]$$
$$=-\mathbb{E}^x\left[\mathbb{E}^{X_{\tau_k+\Delta_k}}\left[\chi_{\{X_{\tau_{k+1}}\leq\bar{x}\}}e^{-\delta\tau_{k+1}}\left(\bar{\Phi}(X_{\tau_{k+1}})-\hat{\bar{\Phi}}(X_{\tau_{k+1}})\right)\right]\right]$$ (71)
$$=-\mathbb{E}^x\left[e^{-\delta(\tau_k+\Delta_k)}\mathbb{E}^{0,X_{\tau_k+\Delta_k}}\left[\chi_{\{X_{\tau_{k+1}-(\tau_k+\Delta_k)}\leq\bar{x}\}}e^{-\delta(\tau_{k+1}-(\tau_k+\Delta_k))}\left(\bar{\Phi}(X_{\tau_{k+1}-(\tau_k+\Delta_k)})-\hat{\bar{\Phi}}(X_{\tau_{k+1}-(\tau_k+\Delta_k)})\right)\right]\right]$$
$$=-\mathbb{E}^x\left[e^{-\delta(\tau_k+\Delta_k)}\mathbb{E}^{0,X_{\tau_k+\Delta_k}}\left[\int_0^{\infty}\lambda e^{-\lambda s}e^{-\delta s}\chi_{\{X_s\leq\bar{x}\}}\left(\bar{\Phi}(X_s)-\hat{\bar{\Phi}}(X_s)\right)\mathrm{d}s\right]\right]$$
$$=-\mathbb{E}^x\left[\lambda e^{-\delta(\tau_k+\Delta_k)}\mathbb{E}^{0,X_{\tau_k+\Delta_k}}\left[\int_0^{+\infty}e^{-(\lambda+\delta)s}\chi_{\{X_s\leq\bar{x}\}}\left(\bar{\Phi}(X_s)-\hat{\bar{\Phi}}(X_s)\right)\mathrm{d}s\right]\right]$$

Summing up last lines of (70) and (71) equals 0, leading to the desired result $\phi(x,\hat{\eta})=\bar{\Phi}(x)$. Because of (68) we get $\Phi=\bar{\Phi}$ in $D$.

□

**(Proof of Proposition 2)**

The FPE seems to have a complicated form, but it is a system of first-order linear differential equations constrained by the condition $p_W, p_N \in C(0,1)$ and the boundary conditions. In what follows, the constants $c_i$ ($i=1,2,3,...$) represent the unknowns to be determined by these constraints. Firstly, (25) and (26) give

$$p_N(x)=c_1 e^{\lambda S^{-1}(x-\bar{x})} \quad \text{and} \quad p_N(x)=\mu S^{-1}c\left(1-\chi_{\{x=1\}}\right)$$ (72)

by the boundary condition (28). The continuity of $p_N$ at $x=\bar{x}$ gives $c_1=\mu S^{-1}c$. We then get

$$p_W(x)=Fc\left(e^{\lambda S^{-1}(x-\bar{x})}-e^{\mu S^{-1}(x-\bar{x})}\right).$$ (73)

By (32) and (33), we get $c_W$ and $c_N$ as the functions of $c$. Finally, the remaining unknown is $c$, which is uniquely determined from the mass conservation condition $\int_D p\mathrm{d}x=1$.

□

**Appendix B: On the simplified model**

We show that there is a unique threshold $\bar{x} \in (0,1)$ if $c, d, \mu^{-1}, \delta > 0$ are small under certain assumption. The system (24) is explicitly written as

$$De^{-\frac{\delta+\lambda}{S}\bar{x}} + A\bar{x} + B + Ce^{-\frac{\delta+\mu}{S}\bar{x}} = \bar{\Phi}(1)e^{\frac{\delta}{S}(1-\bar{x})} \tag{74}$$

and

$$-\frac{\delta+\lambda}{S}De^{-\frac{\delta+\lambda}{S}\bar{x}} + A - \frac{\delta+\mu}{S}Ce^{-\frac{\delta+\mu}{S}\bar{x}} = -\frac{\delta}{S}\bar{\Phi}(1)e^{\frac{\delta}{S}(1-\bar{x})}. \tag{75}$$

All the coefficients including $A$ through $D$ smoothly depend on $c, d, \mu, \delta$. Firstly, we consider a large $\mu$ limit (small delay limit) and then take a small $\delta$ limit (Ergodic limit). Then, we show that the resulting system admits a unique $\bar{x} \in (0,1)$ if $c, d$ are small. By the continuous dependence of the coefficients of (74)-(75) on the parameters $\mu, \delta$, we obtain the unique existence of $\bar{x}$. Notice that the discussion here is rather formal, and does not directly prove the unique existence of the couple $(\bar{\Phi}(1), \bar{x})$.

Here, we use the existing result of the no-delay case [46]. Taking the limit $\mu \to +\infty$ of (74) and (75) with the assumption $\bar{x} \in (0,1)$ yields

$$De^{-\frac{\delta+\lambda}{S}\bar{x}} + A\bar{x} + B = \bar{\Phi}(1)e^{\frac{\delta}{S}(1-\bar{x})} \quad \text{and} \quad -\frac{\delta+\lambda}{S}De^{-\frac{\delta+\lambda}{S}\bar{x}} + A = -\frac{\delta}{S}\bar{\Phi}(1)e^{\frac{\delta}{S}(1-\bar{x})} \tag{76}$$

with

$$A = \frac{1}{\delta+\lambda}\left(-BS + \lambda(\bar{\Phi}(1)+c+d)\right), \quad B = -\frac{\lambda c}{\delta+\lambda}, \quad C = 0, \text{ and } D = \frac{\delta+\lambda-\lambda cS}{(\delta+\lambda)^2}. \tag{77}$$

We now consider the Ergodic limit $\delta \to +0$ of (76) such that (small $\delta$ method [68]):

$$\delta\Phi \to u \quad \text{in } D \text{ as } \delta \to +0. \tag{78}$$

Taking this the limit $\delta \to +0$ in (76) yields

$$\frac{1-cS}{\lambda}e^{-\frac{\lambda}{S}\bar{x}} + c\left(1-\bar{x}+\frac{S}{\lambda}\right) + d = u\left\{\frac{1}{S}(1-\bar{x})+\frac{1}{\lambda}\right\} \quad \text{and} \quad (1-cS)e^{-\frac{\lambda}{S}\bar{x}} + cS = u. \tag{79}$$

Substituting the second equation into the first one yields

$$dS = (u-cS)(1-\bar{x}). \tag{80}$$

By (80), we should have $u > cS$, which is assumed temporally and is justified later. Assuming that $c$ is small such that $cS < 1$, combining (80) with (79) yields

$$(1-\bar{x})e^{-\frac{\lambda}{S}\bar{x}} = \frac{dS}{1-cS}. \tag{81}$$

The left-hand side of (81) is expressed as $F(\bar{x})$ with $F:[0,1] \to \mathbb{R}: F(x) = (1-x)e^{-\frac{\lambda}{S}x}$. Since $F(0) = 1$ and $F(1) = 0$, and $F$ is decreasing, we get the unique existence of $\bar{x} \in (0,1)$ if

$$0 < cS < 1 \quad \text{and} \quad 0 < \frac{dS}{1-cS} < 1, \tag{82}$$

namely if $(c+d)S<1$. The effective Hamiltonian $u$ is then found by substituting this $\bar{x}$ into (80). By $\bar{x} \in (0,1)$, we obtain $u > cS$ under $(c+d)S<1$:

$$u = cS + \frac{dS}{1-\bar{x}} > cS. \tag{83}$$

Consequently, we formally showed the unique existence of $(\bar{\Phi}(1), \bar{x})$ for small $c, d, \mu^{-1}, \delta$. The continuous dependence of the coefficients of (74)-(75) on these parameters justifies the existence of $(\bar{\Phi}(1), \bar{x})$ required in **Proposition 1**.

**Appendix C: Numerical scheme for the FPE**

The numerical scheme for the system (43)-(46) is explained. The spatial discretization of the present scheme employs the cell-vertex discretization [40]. Firstly, the spatial domain $D = [0,1]$ is discretized with the equidistant vertices $\{x_l\}_{l=0,1,2,...,L}$ with $L \in \mathbb{N}$ and $x_l = lL^{-1}$. Assume $L \geq 3$. Set $\Delta x = L^{-1}$ and the cells $\{C_l\}_{l=0,1,2,...,L}$ on which the PDFs are discretized:

$$C_0 = [0, 0.5\Delta x], \quad C_0 = [(l-0.5)\Delta x, (l+0.5)\Delta x] \quad (1 \leq l \leq L-1), \quad C_L = [1-0.5\Delta x, 1]. \tag{84}$$

The lengths of $C_l$ is denoted as $|C_l|$.

The discrete time step is set as $\Delta t > 0$ and the discrete times as $\{t_m\}_{m=0,1,2,...}$ with $t_m = m\Delta t$. Assume that the control is of the form (36). Set the integer-valued variable $\{l_i\}_{i \in M}$ to detect the free boundary $\bar{x}_i \in [0,1]$. More specifically, for each $i \in M$, $l_i$ is the largest integer $l$ such that $l_i \in C_l$. If such an integer $i_l$ is not detected for $i \in M$, then set $l_i = -1$.

The quantity $q$ discretized at $t_m = m\Delta t$ in $C_l$ is represented by super-scripts as $q^{(l,m)}$. For each $m \geq 0$ and $i \in M$, the scheme is presented as follows. The coefficient $c_i = c_i^{(m)}$ is computed as

$$c_i^{(m)} = \sum_{l=0}^{l_i} |C_l| p_{i,W}^{(l,m)}. \tag{85}$$

The fluxes are approximated at each interface between the cells $C_l, C_{l+1}$ in a classical upwind manner:

$$F_{i,N}^{(l,m)} = -S(i,(l+0.5)\Delta x) p_{i,N}^{(l+1,m)} \text{ and } F_{i,W}^{(l,m)} = -S(i,(l+0.5)\Delta x) p_{i,W}^{(l+1,m)} \quad (l=0,1,2,...,L-1). \tag{86}$$

Set the boundary fluxes as

$$F_{i,N}^{(-1,m)} = F_{i,W}^{(-1,m)} = F_{i,N}^{(L,m)} = F_{i,W}^{(L,m)} = 0. \tag{87}$$

The Dirac Delta $\delta_{\{x=1\}}$ is discretized as $\frac{\chi_{\{l=L\}}}{|C_L|}$. For each $i \in M$, the system (43)-(46) is discretized as

$$\frac{p_{i,N}^{(l,m+1)} - p_{i,N}^{(l,m)}}{\Delta t} + \frac{F_{i,N}^{(l,m)} - F_{i,N}^{(l-1,m)}}{|C_l|} + \lambda p_{i,N}^{(l,m)} + \left(\sum_{j \in M, j \neq i} v_{i,j}\right) p_{i,N}^{(l,m)} - \sum_{j \in M, j \neq i} v_{j,i} p_{j,N}^{(l,m)} = 0 \text{ for } 0 \leq l \leq l_i, \tag{88}$$

$$\frac{p_{i,N}^{(l,m+1)} - p_{i,N}^{(l,m)}}{\Delta t} + \frac{F_{i,N}^{(l,m)} - F_{i,N}^{(l-1,m)}}{|C_l|} + \left(\sum_{j \in M, j \neq i} v_{i,j}\right) p_{i,N}^{(l,m)} - \sum_{j \in M, j \neq i} v_{j,i} p_{j,N}^{(l,m)} - \mu c_i^{(m)} \frac{\chi_{\{l=L\}}}{|C_L|} = 0 \quad \text{for} \quad l_i < l \leq L, \quad (89)$$

$$\begin{aligned}
&\frac{p_{i,W}^{(l,m+1)} - p_{i,W}^{(l,m)}}{\Delta t} + \frac{F_{i,W}^{(l,m)} - F_{i,W}^{(l-1,m)}}{|C_l|} \\
&+ \mu p_{i,W}^{(l,m)} - \lambda p_{i,N}^{(l,m)} + \left(\sum_{j \in M, j \neq i} v_{i,j}\right) p_{i,W}^{(l,m)} - \sum_{j \in M, j \neq i} v_{j,i} p_{j,W}^{(l,m)} = 0
\end{aligned} \quad \text{for} \quad 0 \leq l \leq l_i, \quad (90)$$

and

$$p_{i,W}^{(l,m)} = 0 \quad \text{for} \quad l_i < l \leq L \quad \text{and} \quad F_{i,W}^{(l,m)} = 0 \quad \text{for} \quad l_i < l \leq L \quad (91)$$

The boundary condition $p_{i,N}(t,1) = 0$ is not explicitly considered because the numerical solutions are cell-averaged values but not pointwise values. The numerically computed PDFs $p_{i,W}^{(l,m)}$ and $p_{i,N}^{(l,m)}$ are then computed starting from a non-negative initial condition such that

$$\sum_{i \in M} \sum_{l=0}^{L} |C_l| \left( p_{i,N}^{(l,m)} + p_{i,W}^{(l,m)} \right) = 1 \quad \text{for} \quad m = 0 \quad \text{with (91).} \quad (92)$$

The above-presented scheme is possibly diffusive in applications, and is therefore equipped with a WENO reconstruction. More specifically, the fluxes $F_{i,N}^{(l,m)}$ and $F_{i,W}^{(l,m)}$ are evaluated using the reconstructed values based on Falcone and Kalise [44] for $1 \leq l \leq L-2$, while other discretization procedures are unchanged. The fluxes for $l = 0, L-1$ are still the above-presented upwind one because WENO reconstructions are in general not applicable to evaluating fluxes near the boundaries.

Finally, we check that the above-presented scheme is conservative, namely that the scheme satisfies the property (92) for any $m$. We prove the property for $m = 1$ because the proof for $m \geq 2$ is essentially the same. A direct computation shows

$$\begin{aligned}
\sum_{i \in M} \sum_{l=0}^{L} |C_l| \left( p_{i,N}^{(l,1)} + p_{i,W}^{(l,1)} \right) &= \sum_{i \in M} \sum_{l=0}^{l_i} |C_l| p_{i,N}^{(l,1)} + \sum_{i \in M} \sum_{l=l_i+1}^{L} |C_l| p_{i,W}^{(l,1)} + \sum_{i \in M} \sum_{l=0}^{l_i} |C_l| p_{i,N}^{(l,1)} + \sum_{i \in M} \sum_{l=l_i+1}^{L} |C_l| p_{i,W}^{(l,1)} \quad m \in \mathbb{N} \\
&= \sum_{i \in M} \sum_{l=0}^{l_i} |C_l| p_{i,N}^{(l,1)} + \sum_{i \in M} \sum_{l=0}^{l_i} |C_l| p_{i,N}^{(l,1)} + \sum_{i \in M} \sum_{l=l_i+1}^{L} |C_l| p_{i,W}^{(l,1)}
\end{aligned} \quad (93)$$

by (91). The first through the third summations in the right-hand side of (93) are evaluated as

$$\begin{aligned}
\sum_{i \in M} \sum_{l=0}^{l_i} |C_l| p_{i,N}^{(l,1)} &= \sum_{i \in M} \sum_{l=0}^{l_i} |C_l| p_{i,N}^{(l,0)} - \sum_{i \in M} \sum_{l=0}^{l_i} \left( F_{i,N}^{(l,0)} - F_{i,N}^{(l-1,0)} \right) \Delta t \\
&+ \sum_{i \in M} \sum_{l=0}^{l_i} |C_l| \left[ -\lambda p_{i,N}^{(l,0)} \Delta t - \left( \sum_{j \in M, j \neq i} v_{i,j} \right) p_{i,N}^{(l,0)} \Delta t + \sum_{j \in M, j \neq i} v_{j,i} p_{j,N}^{(l,0)} \Delta t \right]
\end{aligned}, \quad (94)$$

$$\begin{aligned}
\sum_{i \in M} \sum_{l=l_i+1}^{L} |C_l| p_{i,W}^{(l,1)} &= \sum_{i \in M} \sum_{l=l_i+1}^{L} |C_l| p_{i,W}^{(l,0)} - \sum_{i \in M} \sum_{l=l_i+1}^{L} \left( F_{i,N}^{(l,0)} - F_{i,N}^{(l-1,0)} \right) \Delta t \\
&+ \sum_{i \in M} \sum_{l=l_i+1}^{L} |C_l| \left[ -\left( \sum_{j \in M, j \neq i} v_{i,j} \right) p_{i,N}^{(l,0)} \Delta t + \sum_{j \in M, j \neq i} v_{j,i} p_{j,N}^{(l,0)} \Delta t + \mu c_i^{(0)} \frac{\chi_{\{l=L\}}}{|C_L|} \Delta t \right]
\end{aligned}, \quad (95)$$

and

$$\sum_{i\in M}\sum_{l=0}^{l_i}|C_l|p_{i,N}^{(l,1)} = \sum_{i\in M}\sum_{l=0}^{l_i}|C_l|p_{i,W}^{(l,0)} - \sum_{i\in M}\sum_{l=0}^{l_i}\left(F_{i,W}^{(l,0)} - F_{i,W}^{(l-1,0)}\right)\Delta t$$
$$+ \sum_{i\in M}\sum_{l=0}^{l_i}|C_l|\left[\lambda p_{i,N}^{(l,0)}\Delta t - \mu p_{i,W}^{(l,0)}\Delta t - \left(\sum_{j\in M, j\neq i}v_{i,j}\right)p_{i,W}^{(l,0)}\Delta t + \sum_{j\in M, j\neq i}v_{j,i}p_{j,W}^{(l,0)}\Delta t\right]. \quad (96)$$

Now, we have

$$-\sum_{i\in M}\sum_{l=0}^{l_i}\left(F_{i,N}^{(l,0)} - F_{i,N}^{(l-1,0)}\right) - \sum_{i\in M}\sum_{l=l_i+1}^{L}\left(F_{i,N}^{(l,0)} - F_{i,N}^{(l-1,0)}\right)\Delta t - \sum_{i\in M}\sum_{l=0}^{l_i}\left(F_{i,W}^{(l,0)} - F_{i,W}^{(l-1,0)}\right)$$
$$= -\sum_{i\in M}\sum_{l=0}^{L}\left(F_{i,N}^{(l,0)} - F_{i,N}^{(l-1,0)}\right) - \sum_{i\in M}\sum_{l=0}^{L}\left(F_{i,W}^{(l,0)} - F_{i,W}^{(l-1,0)}\right) \quad (97)$$
$$= 0$$

by (87) and (91). In addition, we get

$$\sum_{i\in M}\sum_{l=0}^{l_i}|C_l|\left[-\left(\sum_{j\in M, j\neq i}v_{i,j}\right)p_{i,N}^{(l,0)}\Delta t + \sum_{j\in M, j\neq i}v_{j,i}p_{j,N}^{(l,0)}\Delta t\right]$$
$$= \Delta t\sum_{l=0}^{l_i}|C_l|\sum_{i\in M}\left[-\left(\sum_{j\in M, j\neq i}v_{i,j}\right)p_{i,N}^{(l,0)} + \sum_{j\in M, j\neq i}v_{j,i}p_{j,N}^{(l,0)}\right]$$
$$= \Delta t\sum_{l=0}^{l_i}|C_l|\left(-\sum_{i\in M}\sum_{j\in M, j\neq i}v_{i,j}p_{i,N}^{(l,0)} + \sum_{i\in M}\sum_{j\in M, j\neq i}v_{j,i}p_{j,N}^{(l,0)}\right). \quad (98)$$
$$= \Delta t\sum_{l=0}^{l_i}|C_l|\left(-\sum_{i,j\in M, i\neq j}v_{i,j}p_{i,N}^{(l,0)} + \sum_{j,i\in M, j\neq i}v_{j,i}p_{j,N}^{(l,0)}\right)$$
$$= \Delta t\sum_{l=0}^{l_i}|C_l|\left(-\sum_{i,j\in M, i\neq j}v_{i,j}p_{i,N}^{(l,0)} + \sum_{i,j\in M, i\neq j}v_{i,j}p_{i,N}^{(l,0)}\right)$$
$$= 0$$

Similarly, we get

$$\sum_{i\in M}\sum_{l=l_{i+1}}^{L}|C_l|\left[-\left(\sum_{j\in M, j\neq i}v_{i,j}\right)p_{i,N}^{(l,0)}\Delta t + \sum_{j\in M, j\neq i}v_{j,i}p_{j,N}^{(l,0)}\Delta t\right] = 0 \quad (99)$$

and

$$\sum_{i\in M}\sum_{l=0}^{l_i}|C_l|\left[-\left(\sum_{j\in M, j\neq i}v_{i,j}\right)p_{i,W}^{(l,0)}\Delta t + \sum_{j\in M, j\neq i}v_{j,i}p_{j,W}^{(l,0)}\Delta t\right] = 0 \quad (100)$$

Furthermore, we have

$$\sum_{i\in M}\sum_{l=0}^{l_i}|C_l|\mu p_{i,W}^{(l,0)}\Delta t = \Delta t\sum_{i\in M}\mu c_i^{(0)} = \sum_{i\in M}\sum_{l=l_i+1}^{L}|C_l|\mu c_i^{(0)}\frac{\chi_{\{l=L\}}}{|C_L|}\Delta t. \quad (101)$$

We then arrive at the desired result

$$\begin{aligned}
\sum_{i\in M}\sum_{l=0}^{L}|C_l|\left(p_{i,N}^{(l,1)}+p_{i,W}^{(l,1)}\right) &= \sum_{i\in M}\sum_{l=0}^{l_i}|C_l|p_{i,N}^{(l,0)} + \sum_{i\in M}\sum_{l=0}^{l_i}|C_l|p_{i,N}^{(l,0)} + \sum_{i\in M}\sum_{l=l_i+1}^{L}|C_l|p_{i,W}^{(l,0)} \\
&= \sum_{i\in M}\sum_{l=0}^{l_i}|C_l|p_{i,N}^{(l,0)} + \sum_{i\in M}\sum_{l=l_i+1}^{L}|C_l|p_{i,W}^{(l,0)} + \sum_{i\in M}\sum_{l=0}^{l_i}|C_l|p_{i,N}^{(l,0)} + \sum_{i\in M}\sum_{l=l_i+1}^{L}|C_l|p_{i,W}^{(l,0)} \quad , (102)\\
&= \sum_{i\in M}\sum_{l=0}^{L}|C_l|\left(p_{i,N}^{(l,0)}+p_{i,W}^{(l,0)}\right) \\
&= 1
\end{aligned}$$

and obtain the following proposition.

*Proposition C*

*We have the conservation property (92) for $m\in\mathbb{N}$.*

**Acknowledgements**

JSPS Research Grant No. 19H03073 and Kurita Water and Environment Foundation Grant No. 19B018 support this research. We would like to thank the two anonymous reviewers for their careful reading and suggestions that improved the contents of the manuscript.

**References**

[1] Sethi, S. P. (2019). Optimal Control Theory. Springer, Cham.

[2] Øksendal, B., & Sulem, A. (2019). Applied Stochastic Control of Jump Diffusions. Springer, Cham.

[3] Aïd, R., Basei, M., & Pham, H. (2019). A McKean–Vlasov approach to distributed electricity generation development. Mathematical Methods of Operations Research, 1-42. https://doi.org/10.1007/s00186-019-00692-8

[4] Kvamsdal, S. F., Sandal, L. K., & Poudel, D. (2020). Ecosystem wealth in the Barents Sea. Ecological Economics, 171, 106602. https://doi.org/10.1016/j.ecolecon.2020.106602

[5] Hening, A., Tran, K. Q., Phan, T. T., & Yin, G. (2019). Harvesting of interacting stochastic populations. Journal of Mathematical Biology, 79(2), 533-570. https://doi.org/10.1007/s00285-019-01368-x

[6] Zou, X., & Wang, K. (2016). Optimal harvesting for a stochastic Lotka–Volterra predator-prey system with jumps and nonselective harvesting hypothesis. Optimal Control Applications and Methods, 37(4), 641-662. https://doi.org/10.1002/oca.2185

[7] Wang, L., Yuan, J., Wu, C., & Wang, X. (2019). Practical algorithm for stochastic optimal control problem about microbial fermentation in batch culture. Optimization Letters, 13(3), 527-541. https://doi.org/10.1007/s11590-017-1220-z

[8] Yoshioka, H., & Yoshioka, Y. (2019). Modeling stochastic operation of reservoir under ambiguity with an emphasis on river management. Optimal Control Applications and Methods, 40(4), 764-790. https://doi.org/10.1002/oca.2510

[9] Chen, S., Fu, R., Wedge, L., & Zou, Z. (2020). Consumption and portfolio decisions with uncertain lifetimes. Mathematics and Financial Economics, 1-39. https://doi.org/10.1007/s11579-020-00263-0

[10] Han, X., & Liang, Z. (2020). Optimal reinsurance and investment in danger-zone and safe-region. Optimal Control Applications and Methods, 41(3), 765-792. https://doi.org/10.1002/oca.2568

[11] Cao, Y., & Duan, Y. (2020). Joint production and pricing inventory system under stochastic reference price effect. Computers & Industrial Engineering, 106411. https://doi.org/10.1016/j.cie.2020.106411

[12] Ouaret, S., Kenné, J. P., & Gharbi, A. (2019). Production and replacement planning of a deteriorating remanufacturing system in a closed-loop configuration. Journal of Manufacturing Systems, 53, 234-248. https://doi.org/10.1016/j.jmsy.2019.09.006

[13] Øksendal, B. (2003). Stochastic Differential Equations. Springer, Berlin, Heidelberg.

[14] Yin, G. G., & Zhu, C. (2009). Hybrid switching diffusions: properties and applications. Springer, New York, Dordrecht, Heidelberg, London.

[15] Park, S. (2020). Verification theorems for models of optimal consumption and investment with annuitization. Mathematical Social Sciences, 103, 36-44. https://doi.org/10.1016/j.mathsocsci.2019.11.002

[16] Kharroubi, I., Lim, T., & Mastrolia, T. (2020). Regulation of renewable resource exploitation. SIAM Journal on Control and Optimization, 58(1), 551-579. https://doi.org/10.1137/19M1265740

[17] Yaegashi, Y., Yoshioka, H., Unami, K., & Fujihara, M. (2018). A singular stochastic control model for sustainable population management of the fish-eating waterfowl Phalacrocorax carbo. Journal of environmental management, 219, 18-27. https://doi.org/10.1016/j.jenvman.2018.04.099

[18] Xu, L., Xu, S., & Yao, D. (2020). Maximizing expected terminal utility of an insurer with high gain tax by investment and reinsurance. Computers & Mathematics with Applications, 79(3), 716-734. https://doi.org/10.1016/j.camwa.2019.07.023

[19] Forsyth, P., & Labahn, G. (2019). ε-monotone Fourier methods for optimal stochastic control in finance. Journal of Computational Finance, 22(4). https://ssrn.com/abstract=3341491

[20] Picarelli, A., & Reisinger, C. (2020). Probabilistic error analysis for some approximation schemes to optimal control problems. Systems & Control Letters, 137, 104619. https://doi.org/10.1016/j.sysconle.2019.104619

[21] Yoshioka, H., & Tsujimura, M. (2020). Analysis and computation of an optimality equation arising in an impulse control problem with discrete and costly observations. Journal of Computational and Applied Mathematics, 366, 112399. https://doi.org/10.1016/j.cam.2019.112399

[22] Yoshioka, H., Yoshioka, Y., Yaegashi, Y., Tanaka, T., Horinouchi, M., & Aranishi, F. (2020). Analysis and computation of a discrete costly observation model for growth estimation and management of biological resources. Computers & Mathematics with Applications, 79(4), 1072-1093. https://doi.org/10.1016/j.camwa.2019.08.017

[23] Yoshioka, H., Tsujimura, M., Hamagami, K., and Yoshioka Y. (2020). A hybrid stochastic river environmental restoration modeling with discrete and costly observations, Optimal Control Applications and Methods. https://doi.org/10.1002/OCA.2616. in press.

[24] Dyrssen, H., & Ekström, E. (2018). Sequential testing of a Wiener process with costly observations. Sequential Analysis, 37(1), 47-58. https://doi.org/10.1080/07474946.2018.1427973

[25] Winkelmann, S., Schütte, C., & Kleist, M. V. (2014). Markov control processes with rare state observation: theory and application to treatment scheduling in HIV−1. Communications in Mathematical Sciences, 12(5), 859-877. https://dx.doi.org/10.4310/CMS.2014.v12.n5.a4

[26] Wang, H. (2001). Some control problems with random intervention times. Advances in Applied Probability, 33(2), 404-422. https://doi.org/10.1017/S0001867800010867

[27] Pham, H., & Tankov, P. (2008). A model of optimal consumption under liquidity risk with random trading times. Mathematical Finance: An International Journal of Mathematics, Statistics and Financial Economics, 18(4), 613-627. 10.1111/j.1467-9965.2008.00350.x

[28] Federico, S., Gassiat, P., & Gozzi, F. (2017). Impact of time illiquidity in a mixed market without full observation. Mathematical Finance, 27(2), 401-437. https://doi.org/10.1111/mafi.12101

[29] Boyarchenko, S., & Levendorskiĭ, S. (2019). Optimal Stopping Problems in Lévy Models with Random Observations. Acta Applicandae Mathematicae, 163(1), 19-48. https://doi.org/10.1007/s10440-018-0212-z

[30] Liu, J., Yiu, K. F. C., & Bensoussan, A. (2018). Ergodic control for a mean reverting inventory model. Journal of Industrial & Management Optimization, 14(3), 857. https://doi.org/10.3934/jimo.2017079

[31] Yaegashi, Y., Yoshioka, H., Tsugihashi, K., & Fujihara, M. (2019). An exact viscosity solution to a Hamilton–Jacobi–Bellman quasi-variational inequality for animal population management. Journal of Mathematics in Industry, 9(1), 5. https://doi.org/10.1186/s13362-019-0062-y

[32] Bensoussan, A., Long, H., Perera, S., & Sethi, S. (2012). Impulse control with random reaction periods: A central bank intervention problem. Operations Research Letters, 40(6), 425-430. https://doi.org/10.1016/j.orl.2012.06.012

[33] Perera, S., Buckley, W., & Long, H. (2018). Market-reaction-adjusted optimal central bank intervention policy in a forex market with jumps. Annals of Operations Research, 262(1), 213-238. https://doi.org/10.1007/s10479-016-2297-y

[34] Bruder, B., & Pham, H. (2009). Impulse control problem on finite horizon with execution delay. Stochastic Processes and their Applications, 119(5), 1436-1469. https://doi.org/10.1016/j.spa.2008.07.007

[35] Kharroubi, I., Lim, T., & Vath, V. L. (2019). Optimal exploitation of a resource with stochastic population dynamics and delayed renewal. Journal of Mathematical Analysis and Applications, 477(1), 627-656. https://doi.org/10.1016/j.jmaa.2019.04.052

[36] Øksendal, B., & Sulem, A. (2008). Optimal stochastic impulse control with delayed reaction. Applied Mathematics and Optimization, 58(2), 243-255. https://doi.org/10.1007/s00245-007-9034-5

[37] Perera, S., & Long, H. (2017). An approximation scheme for impulse control with random reaction periods. Operations Research Letters, 45(6), 585-591. https://doi.org/10.1016/j.orl.2017.08.014

[38] Yoshioka, H., Yaegashi, Y., Yoshioka, Y., & Hamagami, K. (2019). Hamilton–Jacobi–Bellman quasi-variational inequality arising in an environmental problem and its numerical discretization. Computers & Mathematics with Applications, 77(8), 2182-2206. https://doi.org/10.1016/j.camwa.2018.12.004

[39] Bertucci, C. (2020). Fokker-Planck equations of jumping particles and mean field games of impulse

control. In Annales de l'Institut Henri Poincaré C, Analyse non linéaire. https://doi.org/10.1016/j.anihpc.2020.04.006

[40] Yaegashi, Y., Yoshioka, H., Tsugihashi, K., & Fujihara, M. (2019). Analysis and computation of probability density functions for a 1-D impulsively controlled diffusion process. Comptes Rendus Mathematique, 357(3), 306-315. https://doi.org/10.1016/j.crma.2019.02.007

[41] Evans, M. R., Majumdar, S. N., & Schehr, G. (2020). Stochastic resetting and applications. Journal of Physics A: Mathematical and Theoretical, 53(19), 193001. https://doi.org/10.1088/1751-8121/ab7cfe

[42] Bartlett, M. S., Porporato, A., & Rondoni, L. (2019). Jump processes with deterministic and stochastic controls. Physical Review E, 100(4), 042133. https://doi.org/10.1103/PhysRevE.100.042133

[43] Carlini, E., Ferretti, R., & Russo, G. (2005). A weighted essentially nonoscillatory, large time-step scheme for Hamilton--Jacobi equations. SIAM Journal on Scientific Computing, 27(3), 1071-1091. https://doi.org/10.1137/040608787

[44] Falcone, M., & Kalise, D. (2013). A high-order semi-Lagrangian/finite volume scheme for Hamilton-Jacobi-Isaacs equations. In IFIP Conference on System Modeling and Optimization (pp. 105-117). Springer, Berlin, Heidelberg. https://doi.org/10.1007/978-3-662-45504-3_10

[45] Cortes, J. (2008). Discontinuous dynamical systems. IEEE Control Systems Magazine, 28(3), 36-73. https://doi.org/10.1109/MCS.2008.919306

[46] Yoshioka H. (2020) River environmental restoration based on random observations of a non-smooth stochastic dynamical system. Preprint. https://arxiv.org/submit/3171356

[47] Zhu, C., Yin, G., & Baran, N. A. (2015). Feynman–Kac formulas for regime-switching jump diffusions and their applications. Stochastics An International Journal of Probability and Stochastic Processes, 87(6), 1000-1032. https://doi.org/10.1080/17442508.2015.1019884

[48] Annunziato, M., & Borzì, A. (2014). Optimal control of a class of piecewise deterministic processes. European Journal of Applied Mathematics, 25(1), 1-25. https://doi.org/10.1017/S0956792513000259

[49] Ramponi, A. (2011). Mixture dynamics and regime switching diffusions with application to option pricing. Methodology and Computing in Applied Probability, 13(2), 349-368. https://doi.org/10.1007/s11009-009-9155-1

[50] Jiang, G. S., & Peng, D. (2000). Weighted ENO schemes for Hamilton--Jacobi equations. SIAM Journal on Scientific Computing, 21(6), 2126-2143. https://doi.org/10.1137/S106482759732455X

[51] Shu, C. W. (2009). High order weighted essentially nonoscillatory schemes for convection dominated problems. SIAM Review, 51(1), 82-126. https://doi.org/10.1137/070679065

[52] Oberman, A. M. (2006). Convergent difference schemes for degenerate elliptic and parabolic equations: Hamilton--Jacobi equations and free boundary problems. SIAM Journal on Numerical Analysis, 44(2), 879-895. https://doi.org/10.1137/S0036142903435235

[53] Xu, S., Chen, M., Liu, C., Zhang, R., & Yue, X. (2019). Behavior of different numerical schemes for random genetic drift. BIT Numerical Mathematics, 59(3), 797-821. https://doi.org/10.1007/s10543-019-00749-4

[54] Zhang, X. F., Yan, H. C., Yue, Y., & Xu, Q. X. (2019). Quantifying natural and anthropogenic impacts

on runoff and sediment load: An investigation on the middle and lower reaches of the Jinsha River Basin. Journal of Hydrology: Regional Studies, 25, 100617. https://doi.org/10.1016/j.ejrh.2019.100617

[55] Katz, S. B., Segura, C., & Warren, D. R. (2018). The influence of channel bed disturbance on benthic Chlorophyll a: A high resolution perspective. Geomorphology, 305, 141-153. https://doi.org/10.1016/j.geomorph.2017.11.010

[56] Neverman, A. J., Death, R. G., Fuller, I. C., Singh, R., & Procter, J. N. (2018). Towards mechanistic hydrological limits: a literature synthesis to improve the study of direct linkages between sediment transport and periphyton accrual in gravel-bed rivers. Environmental Management, 62(4), 740-755. https://doi.org/10.1007/s00267-018-1070-1

[57] Brousse, G., Arnaud-Fassetta, G., Liébault, F., Bertrand, M., Melun, G., Loire, R., ... & Borgniet, L. (2019). Channel response to sediment replenishment in a large gravel-bed river: The case of the Saint-Sauveur dam in the Buëch River (Southern Alps, France). River Research and Applications. https://doi.org/10.1002/rra.3527

[58] Stähly, S., Franca, M. J., Robinson, C. T., & Schleiss, A. J. (2019). Sediment replenishment combined with an artificial flood improves river habitats downstream of a dam. Scientific Reports, 9(1), 1-8. https://doi.org/10.1038/s41598-019-41575-6

[59] Szymkiewicz, R. (2010). Numerical Modeling in Open Channel Hydraulics. Springer, Dordrecht.

[60] Meyer-Peter, E., & Müller, R. (1948). Formulas for bed-load transport. In IAHSR 2nd meeting, Stockholm, appendix 2. IAHR.

[61] Wong, M., & Parker, G. (2006). Reanalysis and correction of bed-load relation of Meyer-Peter and Müller using their own database. Journal of Hydraulic Engineering, 132(11), 1159-1168. https://doi.org/10.1061/(ASCE)0733-9429(2006)132:11(1159)

[62] Turner, S. W. D., & Galelli, S. (2016). Regime-shifting streamflow processes: Implications for water supply reservoir operations. Water Resources Research, 52(5), 3984-4002. https://doi.org/10.1002/2015WR017913

[63] Crandall, M. G., Ishii, H., & Lions, P. L. (1992). User's guide to viscosity solutions of second order partial differential equations. Bulletin of the American Mathematical Society, 27(1), 1-67. https://doi.org/10.1090/S0273-0979-1992-00266-5

[64] Barles, G., & Souganidis, P. E. (1991). Convergence of approximation schemes for fully nonlinear second order equations. Asymptotic Analysis, 4(3), 271-283. https://doi.org/10.3233/ASY-1991-4305

[65] do Val, J. B. R., Guillotreau, P., & Vallée, T. (2019). Fishery management under poorly known dynamics. European Journal of Operational Research, 279(1), 242-257. https://doi.org/10.1016/j.ejor.2019.05.016

[66] Li, W., Huang, L., & Wang, J. (2020). Dynamic analysis of discontinuous plant disease models with a non-smooth separation line. Nonlinear Dynamics, 99(2), 1675-1697. https://doi.org/10.1007/s11071-019-05384-w

[67] Yu, Z. (2012). The stochastic maximum principle for optimal control problems of delay systems

involving continuous and impulse controls. Automatica, 48(10), 2420-2432. https://doi.org/10.1016/j.automatica.2012.06.082

[68] Qian, J. L. (2003). Two Approximations for Effective Hamiltonians Arising from Homogenization of Hamilton-Jacobi Equations. UCLA CAM report. Los Angeles, CA: University of California; 2003.